\documentclass{statsoc}
\usepackage{amssymb,amsmath,ifthen,bbm,subfigure,verbatim}
\usepackage{algorithm}
\usepackage{algorithmicx}
\usepackage{algpseudocode}
\usepackage{latexsym}
\usepackage{natbib}

\usepackage{graphics}
\usepackage[dvips]{epsfig}

\numberwithin{equation}{section}

\newcommand{\APF}{APF\ }
\newcommand{\bddonX}[1]{\mathcal{B}_{\mathrm{b}}(\Xset^{#1})}
\newcommand{\define}{\triangleq}
\newcommand{\E}{\mathbb{E}}
\newcommand{\eqsp}{\;}
\newcommand{\filt}[1]{\mathcal{F}_{#1}}
\newcommand{\filthat}[1]{\hat{\mathcal{F}}_{#1}}
\newcommand{\filttilde}[1]{\tilde{\mathcal{F}}_{#1}}
\newcommand{\fop}[2]{\Psi_{#1:#2}}
\newcommand{\fstimp}[1]{\varphi_{#1}^N}
\newcommand{\fstimpfunc}[1]{\fstimpfuncletter_{#1}}
\newcommand{\fstimpfuncletter}{t}
\newcommand{\fsttarg}[1]{\alpha_{#1}^N}
\newcommand{\fstwgt}[2]{{\tau}_{#1}^{N,#2}}
\newcommand{\hd}[1]{q_{#1}}
\newcommand{\hk}[1]{Q_{#1}}
\newcommand{\iid}{i.i.d.}
\newcommand{\ind}{\mathbbm{1}}
\newcommand{\lb}{\epsilon_-}
\newcommand{\lf}[1]{g_{#1}}
\newcommand{\Lp}[1]{\mathsf{L}^{#1}}
\newcommand{\meanopt}[1]{\tilde{m}_{#1}}
\newcommand{\meanoptSV}[1]{\bar{m}_{#1}}
\newcommand{\opt}[1]{H_{#1}}
\newcommand{\osc}{\operatorname{osc}}
\newcommand{\overpil}[1]
  {\raisebox{1.5ex}{ $\underrightarrow{\text{\ \scriptsize #1\ }}$ }}
\newcommand{\Prob}{\operatorname{\mathbb{P}}}
\newcommand{\proj}[1]{\pi_{#1}}
\newcommand{\prop}[1]{R_{#1}}
\newcommand{\propdens}[1]{r_{#1}}
\newcommand{\parti}[3]
{\ifthenelse{\equal{#3}{}}{\xi_{#1}^{N,#2}}{\xi_{#1}^{N,#2}(#3)}}
\newcommand{\parthat}[3]
{\ifthenelse{\equal{#3}{}}{\hat{\xi}_{#1}^{N,#2}}{\hat{\xi}_{#1}^{N,#2}(#3)}}
\newcommand{\partsmooth}[2]{\phi^N_{#1}}
\newcommand{\partsmoothtilde}[2]{\tilde{\phi}^N_{#1}}
\newcommand{\pathprop}[1]{R_{#1}^{\mathrm{p}}}
\newcommand{\parttilde}[3]
{\ifthenelse{\equal{#3}{}}{\tilde{\xi}_{#1}^{N,#2}}{\tilde{\xi}_{#1}^{N,#2}(#3)}}
\newcommand{\R}{\mathbb{R}}
\newcommand{\refm}{\mu}
\newcommand{\PSimpfunct}[1]{\fstimpfunc{#1}^{\mathrm{P\&S}}}
\newcommand{\sdimpfunc}[1]{\sdimpfuncletter_{#1}}

\newcommand{\sdtarg}[1]{\beta_{#1}^N}
\newcommand{\sigmaopt}[1]{\tilde{\sigma}_{#1}}
\newcommand{\sigmaoptSV}[1]{\bar{\sigma}_{#1}}
\newcommand{\smooth}[2]{\phi_{#1}}
\newcommand{\smoothmixture}[2]{\bar{\phi}^N_{#1}}
\newcommand{\smoothop}[1]{\Phi_{#1}}
\newcommand{\supnm}[2]{\left \| #2 \right \|_{#1, \infty}}
\newcommand{\SSAPF}{SSAPF\ }
\newcommand{\TSS}{TSSPF\ }
\newcommand{\ub}{\epsilon_+}
\newcommand{\ud}{\mathrm{d}}
\newcommand{\Xset}{\mathsf{X}}
\newcommand{\Xsigm}{\mathcal{X}}
\newcommand{\uk}[1]{H^{\mathrm{u}}_{#1}}
\newcommand{\vect}[1]{#1}
\newcommand{\yrefm}{\lambda}
\newcommand{\Yset}{\mathsf{Y}}
\newcommand{\Ysigm}{\mathcal{Y}}
\newcommand{\wgt}[2]{\omega_{#1}^{N,#2}}
\newcommand{\wgtsum}[1]{\Omega_{#1}^N}
\newcommand{\wgtsumtilde}[1]{\tilde{\Omega}_{#1}^N}
\newcommand{\wgttilde}[2]{\tilde{\omega}_{#1}^{N,#2}}
\newcommand{\wrt}{with respect to\ }

\newcounter{hyp}
\newenvironment{hyp}[1]{\refstepcounter{hyp}\it\begin{itemize}\item[{\bf
      (A\arabic{hyp})}] \label{hyp:#1}}{\end{itemize}}
\newcommand{\refhyp}[1]{{\bf (A\ref{hyp:#1})}}
\renewenvironment{proof}[2][Proof]
  {\vspace{10pt plus 3pt minus 2pt}
  {\em #1{}.{ }}#2}
  {\hfill$\Box$\vspace{10pt plus 3pt minus 2pt}\\}

\newtheorem{corollary}{Corollary}[section]
\newtheorem{definition}{Definition}[section]
\newtheorem{lemma}{Lemma}[section]
\newtheorem{theorem}{Theorem}[section]

\title[On the auxiliary particle filter]{On the auxiliary particle filter}

\author[Douc {\it et al.}]{R. Douc}
\address{Ecole Polytechnique,
         Paris,
         France.}
\email{}

\author[Douc {\it et al.}]{\'{E}. Moulines and J. Olsson}
\address{Ecole Nationale Sup\'{e}rieure des T\'{e}l\'{e}communications,
         Paris,
         France.}
\email{\{moulines, olsson\}@enst.tsi.fr}

\begin{document}

 \begin{abstract}
 In this article we study asymptotic properties of weighted samples
 produced by the auxiliary particle filter (APF) proposed by
 \citet{pitt:shephard:1999}. Besides establishing a central limit
 theorem (CLT) for smoothed particle estimates, we also derive bounds on the
 $\Lp{p}$ error and bias of the same for a finite particle sample
 size. By examining the recursive formula for the asymptotic variance
 of the CLT we identify first-stage importance weights for which the
 increase of asymptotic variance at a single iteration of the algorithm
 is minimal. In the light of these findings, we discuss and
 demonstrate on several examples how the \APF algorithm can be improved. 
 \end{abstract}

 \section{Introduction} \label{section:introduction}
 \setcounter{equation}{0}
 In this paper we consider a \emph{state space model} where a sequence
 $Y \define \{ Y_k \}_{k = 0}^\infty$ is modeled as a 
 noisy observation of a Markov chain $X \define \{ X_k \}_{k =
   0}^\infty$, called the \emph{state sequence}, which is hidden. The
 observed values of $Y$ are conditionally independent given the hidden
 states $X$ and the corresponding conditional distribution of $Y_k$
 depends on $X_k$ only. When operating on a model of this form the
 \emph{joint smoothing distribution}, that is, the joint 
distribution of $(X_0, \ldots, X_n)$ given $(Y_0, \ldots, Y_n)$, and
its marginals will be of interest. Of particular interest is the
\emph{filter distribution}, defined as the marginal of this law \wrt
the component $X_n$ is referred to. Computing these 
 posterior distributions will be the key issue when filtering the hidden
 states as well as performing inference on unknown model
 parameters. The posterior distribution can be recursively updated as new
observations become available---making single-sweep processing of the
data possible---by means of the so-called \emph{smoothing
  recursion}. However, in general this recursion cannot be applied
directly since it involves the evaluation of complicated high-dimensional
integrals. In fact, closed form solutions are obtainable only for
linear/Gaussian models (where the solutions are acquired using
the \emph{disturbance smoother}) and models where the state
space of the latent Markov chain is finite.  

\emph{Sequential Monte Carlo} (SMC) \emph{methods}, often alternatively
termed \emph{particle filters}, provide a helpful tool for computing
approximate solutions to the smoothing recursion for general state
space models, and the field has seen a drastic increase in interest
over recent years. These methods 
are based on the principle of, recursively in time, approximating the
smoothing distribution with the empirical measure
associated with a weighted sample of \emph{particles}. At present  
time there are various techniques for producing and updating such a
particle sample
\citep[see][]{fearnhead:1998,doucet:defreitas:gordon:2001,liu:2001}.  
For a comprehensive treatment of the theoretical aspects of SMC
methods we refer to the work by \citet{delmoral:2004}.

In this article we analyse the \emph{auxiliary particle filter} (APF)
proposed by \citet{pitt:shephard:1999}, which has proved to be one of
the most useful and widely adopted implementations of the SMC
methodology. Unlike the traditional \emph{bootstrap particle filter}
\citep{gordon:salmond:smith:1993}, the APF enables the user to affect
the particle sample allocation by designing freely a set of 
\emph{first-stage importance weights} involved in the selection
procedure. Prevalently, this has been used for assigning large weight
to particles whose offsprings are 
likely to land up in zones of the state space having high posterior
probability. Despite its obvious appeal, it is however not clear how
to optimally exploit this additional degree of freedom.  

In order to better understand this issue, we present an asymptotical
analysis (being a continuation of \citep{douc:moulines:olsson:2006} and based
on recent results by
\citep{chopin:2004,kunsch:2005,douc:moulines:2005} on weighted systems
of particles) of the algorithm. More specifically, we establish CLTs
(Theorems~\ref{th:CLT:TSS} and \ref{th:CLT:SSAPF}), with explicit
expressions of the asymptotic variances, for two different versions
(differentiated by the absence/presence of a concluding resampling
pass at the end of each loop) of the algorithm  
under general model specifications. The convergence bear upon an
increasing number of particles, and a recent result in the same
spirit has, independently of \citep{douc:moulines:olsson:2006}, been
stated in the manuscript \citep{doucet:johansen:2007}. Using these
results, we also---and this is the main contribution 
of the paper---identify first-stage importance weights which are
asymptotically most efficient. This result provides important insights
in optimal sample allocation for particle filters in general, and we
also give an interpretation of the finding in terms of variance
reduction for stratified sampling.  

In addition, we prove (utilising a decomposition of
the Monte Carlo error proposed by \citet{delmoral:2004} and refined
by \citet{cappe:douc:moulines:olsson:2006}) time uniform convergence
in $\Lp{p}$ (Theorem~\ref{th:main:deviation:theorem}) under more
stringent assumptions of ergodicity of the conditional hidden
chain. With support of this stability result and the asymptotic
analysis we conclude that inserting a final selection step at the end
of each loop is---at least as long as the number of particles used in
the two stages agree---superfluous, since such an operation
exclusively increases the asymptotic variance. 

Finally, in the implementation section
(Section~\ref{section:implementations}) several heuristics, derived
from the obtained results, for designing efficient first-stage weights
are discussed, and the improvement implied by approximating the asymptotically
optimal first-stage weights is demonstrated on several examples.   

\section{Notation and basic concepts}
\subsection{Model description}
We denote by $(\Xset, \Xsigm)$, $\hk{}$, and $\nu$ the state space,
transition kernel, and initial distribution of $X$, respectively, and
assume that all random variables are defined on a common probability
space $(\Omega, \Prob, \mathcal{A})$. In addition we denote by
$(\Yset, \Ysigm)$ the state space of $Y$ and suppose that there exists
 a measure $\yrefm$ and, for all $x \in \Xset$, a 
 non-negative function $y \mapsto g(y | x)$ such that, for $k \geq 0$,
 $\Prob \left( Y_k \in A | X_k = x \right) = 
 \int_A g(y | x) \, \yrefm(\ud y)$, $A \in \Ysigm$. Introduce, for $i
 \leq j$, the vector notation $\vect{X}_{i:j} \define 
(X_i, \ldots, X_j)$; similar notation will be used for other
quantities. The joint smoothing distribution of denoted by
 \[
 \smooth{n}{n}(A) \define
 \Prob \left(\left. \vect{X}_{0:n} \in A
   \right| \vect{Y}_{0:n} = \vect{y}_{0:n} \right) \eqsp, \quad A \in
 \Xsigm^{\otimes(n+1)} \eqsp,
 \] 
and a straightforward
 application of Bayes's formula shows that 
\begin{equation} \label{eq:smoothing:recursion}
\smooth{k+1}{k+1}(A) = \frac{ \int_A
  \lf{}(y_{k+1} | x_{k+1}) \, \hk{}(x_k, 
  \ud x_{k+1}) \, \smooth{k}{k}(\ud \vect{x}_{0:k}) }{
  \int_{\Xset^{k+2}} \lf{}(y_{k+1} | x'_{k+1}) \, \hk{}(x'_k, \ud 
  x'_{k+1}) \, \smooth{k}{k}(\ud \vect{x}'_{0:k}) } \eqsp,
\end{equation}
for sets $A \in \Xsigm^{\otimes (k+2)}$.
We will throughout this paper assume that we are given a
 sequence $\{y_k ; k \geq 0\}$ of \emph{fixed} observations, and
 write, for $x \in \Xset$, $\lf{k}(x) \define \lf{}(y_k |
 x)$. Moreover, from now on we let the dependence on
 these observations of all other quantities be implicit, and denote, since
 the coming analysis is made exclusively 
 \emph{conditionally} on the given observed record, by $\Prob$
 and $\E$ the conditional probability measure and 
 expectation \wrt these observations. 

\subsection{The auxiliary particle filter}
\setcounter{equation}{0}
\label{section:two-stage:sampling}
Let us recall the \APF algorithm by \citet{pitt:shephard:1999}. Assume
that we at time $k$ have a particle sample $\{ 
(\parti{0:k}{i}{}, \wgt{k}{i}) \}_{i = 1}^N$ (each random variable
$\parti{0:k}{i}{}$ taking values in $\Xset^{k+1}$) 
providing an approximation $\sum_{i=1}^N
\wgt{k}{i} \delta_{\parti{0:k}{i}{}} / \wgtsum{k}$ of the joint
smoothing distribution $\smooth{k}{k}$,
where $\Omega_k^N \define \sum_{i=1}^N 
 \omega_k^{N,i}$ and $\omega_k^{N,i} \geq 0$, $1 \leq i \leq N$. Then,
 when the observation $y_{k+1}$ becomes available, an 
approximation of $\smooth{k+1}{k+1}$ is obtained by plugging
the empirical measure $\partsmooth{k}{k}$ into the recursion
\eqref{eq:smoothing:recursion}, yielding, for $A \in \Xsigm^{\otimes
  (k+1)}$,
\begin{equation*} \label{eq:final:est}
\smoothmixture{k+1}{k+1}(A) \define
\sum_{i=1}^N \frac{\wgt{k}{i} \uk{k}(\parti{0:k}{i}{},
  \Xset^{k+2})}{\sum_{j=1}^N \wgt{k}{j} \uk{k}(\parti{0:k}{j}{},
  \Xset^{k+2})} \opt{k} (\parti{0:k}{i}{}, A) \eqsp, \quad A \in
 \Xsigm^{\otimes(n+1)} \eqsp.
\end{equation*}
Here we have introduced, for $\vect{x}_{0:k}
\in \Xset^{k+1}$ and $A \in \Xsigm^{\otimes (k+1)}$, the unnormalised kernels
 \[
 \uk{k}(\vect{x}_{0:k}, A) \define \int_A \lf{k+1}(x'_{k+1}) \,
 \delta_{\vect{x}_{0:k}}(\ud \vect{x}'_{0:k}) \, \hk{}(x'_k, \ud
 x'_{k+1})
 \]
and $\opt{k}(\vect{x}_{0:k}, A) \define \uk{k}(\vect{x}_{0:k},
  A) / \uk{k}(\vect{x}_{0:k}, \Xset^{k+2})$. Simulating from
  $\opt{k}(\vect{x}_{0:k}, A)$ consists in extending the trajectory
  $\vect{x}_{0:k} \in \Xset^{k+1}$ with an additional component being
  distributed according to the \emph{optimal kernel}, that is, the
  distribution of $X_{k+1}$ conditional on $X_k = x_k$ \emph{and} the
  observation $Y_{k+1} = y_{k+1}$. Now, since we want to form a new
  weighted sample approximating 
$\smooth{k+1}{k+1}$, we need to find a convenient mechanism for
sampling from $\smoothmixture{k+1}{k+1}$ given $\{
(\parti{0:k}{i}{}, \wgt{k}{i}) \}_ {i = 1}^N$.
In most cases cases it is possible---but generally
computationally expensive---to simulate from
$\smoothmixture{k+1}{k+1}$ directly using
\emph{auxiliary accept-reject sampling}
\citep[see][]{hurzeler:kunsch:1998,kunsch:2005}. A computationally
cheaper \citep[see][p.~1988, for a discussion of the acceptance
probability associated with the auxiliary accept-reject sampling
approach]{kunsch:2005} solution consists in producing a weighted
sample approximating $\smoothmixture{k+1}{k+1}$ by sampling from the
importance sampling distribution
\begin{equation*}  \label{eq:imp:sampling:distr}
\rho_{k+1}^N (A) \define \sum_{i=1}^N
 \frac{\wgt{k}{i} \fstwgt{k}{i}}{\sum_{j=1}^N \wgt{k}{j}
   \fstwgt{k}{j}} \pathprop{k} (\parti{0:k}{i}{}, A) \eqsp, \quad A
 \in \Xsigm^{\otimes (k+2)} \eqsp.
 \end{equation*}
 Here $\fstwgt{k}{i}$, $1 \leq i \leq N$, are positive numbers
 referred to as \emph{first-stage weights} \citep[][use the term
 \emph{adjustment multiplier weights}]{pitt:shephard:1999} and
in this article we consider first-stage weights of type
\begin{equation} \label{eq:def:fstimpfunc}
\fstwgt{k}{i} = \fstimpfunc{k} (\parti{0:k}{i}{})
\end{equation}
for some function $\fstimpfunc{k} : \Xset^{k+1} \rightarrow
\R^{+}$. Moreover, the pathwise proposal kernel $\pathprop{k}$ is,
for $\vect{x}_{0:k} \in \Xset^{k+1}$ and $A \in
\Xsigm^{\otimes(k+2)}$, of form
\[
\pathprop{k}(\vect{x}_{0:k},A) = \int_A
\delta_{\vect{x}_{0:k}}(\ud \vect{x}'_{0:k}) \,
\prop{k}(x'_k, \ud x'_{k+1})
\]
with $\prop{k}$ being such that $\hk{}(x, \cdot) \ll \prop{k}(x,
\cdot)$ for all $x \in \Xset$. Thus, a draw from
$\pathprop{k}(\vect{x}_{0:k},\cdot)$ is produced by extending the
trajectory $\vect{x}_{0:k} \in \Xset^{k+1}$ with an additional
component obtained by simulating from $\prop{k}(x_k, \cdot)$.
 It is easily checked that for $\vect{x}_{0:k+1} \in \Xset^{k+2}$,
\begin{equation} \label{eq:def:sdimpfunc}
 \frac{\ud
  \smoothmixture{k+1}{k+1}}{\ud
  \rho_{k+1}^N}(\vect{x}_{0:k+1}) \propto
\sdimpfunc{k+1}(\vect{x}_{0:k+1}) \define \sum_{i=1}^N 
\ind_{\parti{0:k}{i}{}}(\vect{x}_{0:k})
\frac{\lf{k+1}(x_{k+1})}{\fstwgt{k}{i}} \frac{\ud \hk{}(x_k,
  \cdot)}{\ud \prop{k}(x_k, \cdot)}(x_{k+1}) \eqsp.
\end{equation}
An updated weighted particle sample $\{ (\parttilde{0:k+1}{i}{},
\wgttilde{k+1}{i}{}) \}_ {i = 1}^{M_N}$ targeting
$\smoothmixture{k+1}{k+1}$ is hence generated by simulating $M_N$
particles $\parttilde{0:k+1}{i}{}$, $1 \leq i \leq M_N$, from the
proposal $\rho_{k+1}^N$ and associating with these particles the
\emph{second-stage weights} $\wgttilde{k+1}{i} \define
\sdimpfunc{k+1}(\parttilde{0:k+1}{i}{})$, $1 \leq i \leq M_N$. By the
identity function in \eqref{eq:def:sdimpfunc}, only a single term of
the sum will contribute to the second-stage weight of a particle. 

Finally, in an \emph{optional} second-stage resampling pass a
uniformly weighted particle sample $\{ (\parttilde{0:k+1}{i}{}, 1) \}_
{i = 1}^N$, still targeting $\smoothmixture{k+1}{k+1}$, is
obtained by resampling $N$ of the particles $\parttilde{0:k+1}{i}{}$, $1
\leq i \leq M_N$, according to the normalised second-stage
weights. Note that the number of particles in the last two samples,
$M_N$ and $N$, may be different. The procedure is now repeated
recursively (with $\wgt{k+1}{i} \equiv 1$, $1 \leq i \leq N$) and is
initialised by drawing $\parti{0}{i}{}$, $1 \leq i \leq N$, independently from
$\varsigma$, where $\nu \ll \varsigma$, yielding
$\wgt{0}{i} = \sdimpfunc{0}(\parti{0}{i}{})$ with
$\sdimpfunc{0}(x) \define \lf{0}(x) \, \ud \nu / \ud \varsigma (x)$, $x
\in \Xset$. To summarise, we obtain, depending on whether second-stage
resampling is performed or not, the procedures described in
Algorithms~1 and~2.

\begin{algorithm}[h]
\caption{Two-Stage Sampling Particle Filter (TSSPF)}
\label{alg:TSS}
\begin{algorithmic}[1]
\Ensure $\{ (\parti{0:k}{i}{}, \wgt{k}{i}) \}_{i=1
 }^N$ approximates $\smooth{k}{k}$.
 \For{$i = 1, \ldots, M_N$} \Comment{First stage}
 \State draw indices $I_k^{N,i}$ from the set $\{1, \ldots, N
 \}$ multinomially with respect to the normalised
 weights $\wgt{k}{j} \fstwgt{k}{j} / \sum_{\ell=1}^N
 \wgt{k}{\ell} \fstwgt{k}{\ell}$, $1 \leq j \leq N$;
 \State simulate $\parttilde{0:k+1}{i}{k+1} \sim
 \prop{k}[\parti{0:k}{I_k^{N,i}}{k}, \cdot]$, and
 \State set $\parttilde{0:k+1}{i}{} \define [\parti{0:k}{I_k^{N,i}}{},
 \parttilde{0:k+1}{i}{k+1}]$ and $\wgttilde{k+1}{i} \define
 \sdimpfunc{k+1}(\parttilde{0:k+1}{i}{})$.
 \EndFor
 \For{$i = 1, \ldots, N$} \Comment{Second stage}
 \State draw indices $J_{k+1}^{N,i}$ from the set $\{1, \ldots, M_N
 \}$ multinomially with respect to the
 normalised
 weights $\wgttilde{k+1}{j} / \sum_{\ell=1}^N \wgttilde{k+1}{\ell}$, $1 \leq j
 \leq N$, and
 \State set $\parti{0:k+1}{i}{} \define \parttilde{0:k+1}{J_{k+1}^{N,i}}{}$.
 \State Finally, reset the weights: $\wgt{k+1}{i} = 1$.
 \EndFor
 \State Take $\{ (\parti{0:k+1}{i}{}, 1) \}_{i=1
 }^N$ as an approximation of $\smooth{k+1}{k+1}$.
\end{algorithmic}
\end{algorithm}

\begin{algorithm}[h]
\caption{Single-Stage Auxiliary Particle Filter (SSAPF)}
\label{alg:single:step:aux}
\begin{algorithmic}[1]
\Ensure $\{ (\parti{0:k}{i}{}, \wgt{k}{i}) \}_{i=1
}^N$ approximates $\smooth{k}{k}$.
\For{$i = 1, \ldots, N$}
\State draw indices $I_k^{N,i}$ from the set $\{1, \ldots, N
 \}$ multinomially with respect to the normalised
weights $\wgt{k}{j} \fstwgt{k}{j} / \sum_{\ell=1}^N
\wgt{k}{\ell} \fstwgt{k}{\ell}$, $1 \leq j \leq N$;
\State simulate $\parttilde{0:k+1}{i}{k+1} \sim
\prop{k}[\parti{0:k}{I_k^{N,i}}{k}, \cdot]$, and
\State set
$\parttilde{0:k+1}{i}{} \define [\parti{0:k}{I_k^{N,i}}{},
\parttilde{0:k+1}{i}{k+1}]$ and $\wgttilde{k+1}{i} \define
\sdimpfunc{k+1}(\parttilde{0:k+1}{i}{})$.
\EndFor
\State Take $\{ (\parttilde{0:k+1}{i}{}, \wgttilde{k+1}{i}) \}_{i=1
}^N$ as an approximation of $\smooth{k+1}{k+1}$.
\end{algorithmic}
\end{algorithm}

We will use the term \APF as a family name for both these algorithms and refer
to them separately as \emph{two-stage sampling particle filter} (TSSPF)
and \emph{single-stage auxiliary particle filter} (SSAPF).
Note that we by letting $\fstwgt{k}{i} \equiv 1$, $1 \leq i \leq N$, in
Algorithm~\ref{alg:single:step:aux} obtain the bootstrap particle
filter suggested by \citet{gordon:salmond:smith:1993}.

The resampling steps of the \APF can of course be
implemented using techniques (e.g., \emph{residual} or
\emph{systematic} resampling) different from multinomial resampling,
leading to straightforward adaptions not discussed here. We believe
however that the results of the coming analysis are generally
applicable and extendable to a large class of selection schemes.

The issue whether second-stage resampling should be performed or not
has been treated by several authors, and the theoretical results on the
particle approximation stability and asymptotic variance presented in
the next section will indicate that the second-stage selection pass
should, at least for the case $M_N = N$, be canceled, since this
exclusively increases the sampling variance. Thus, the idea
that the second-stage resampling pass is necessary for preventing the
particle approximation from degenerating does not apparently
hold. Recently, a similar conclusion was reached in the manuscript
\citep{doucet:johansen:2007}.

The advantages of the \APF not possessed by
standard SMC methods is the possibility of, firstly, choosing the
first-stage weights $\fstwgt{k}{i}$ arbitrarily and, secondly, letting $N$
and $M_N$ be different (\TSS only). Appealing to common sense,
SMC methods work efficiently when the particle weights are
well-balanced, and \cite{pitt:shephard:1999} propose several
strategies for achieving this by adapting the first-stage weights. In
some cases it is possible to fully adapt the filter to the model (see
Section~\ref{section:implementations}), providing exactly equal
importance weights; otherwise, \cite{pitt:shephard:1999} suggest, in
the case $\prop{k} \equiv \hk{}$ and $\Xset = \R^d$, the generic
first-stage importance weight function $\PSimpfunct{k}(\vect{x}_{0:k})
\define \lf{k+1}[ \int_{\R^d} x' \, \hk{}(x_k, \ud x')]$,
$\vect{x}_{0:k} \in \R^{k+1}$. The analysis that follows will however
show that this way of adapting the first-stage weights is not
necessarily good in terms of asymptotic (as $N$ tends to infinity) 
sample variance; indeed, using first-stage weights given by
$\PSimpfunct{k}$ can be even detrimental for some models.

\section{Bounds and asymptotics for produced approximations}
\label{section:asymptotic:properties}
\setcounter{equation}{0}
\subsection{Asymptotic properties.}
Introduce, for any probability
measure $\mu$ on some measurable space $(\mathsf{E}, \mathcal{E})$ and
$\mu$-measurable function $f$ satisfying $\int_\mathsf{E} |f(x)| \,
\mu(\ud x) < \infty$, the notation $\mu f \define \int_\mathsf{E} f(x) \,
\mu(\ud x)$. Moreover, for any two transition kernels $K$ and $T$ from
$(\mathsf{E}_1, \mathcal{E}_1)$ to $(\mathsf{E}_2, \mathcal{E}_2)$ and
$(\mathsf{E}_2, \mathcal{E}_2)$ to $(\mathsf{E}_3, \mathcal{E}_3)$,
respectively, we define the product transition kernel $KT(x, A)
\define \int_{\mathsf{E}_2}  T(z, A)\, K(x, \ud z)$, for $x \in
\mathsf{E}_1$ and $A \in \mathcal{E}_3$. A set $\mathsf{C}$ of
real-valued functions on $\Xset^m$ is said to be \emph{proper} if the
following conditions hold: $\textbf{i)}$ $\mathsf{C}$ is a linear
space; $\textbf{ii)}$ if $g \in \mathsf{C}$ and $f$ is measurable with
$|f| \leq |g|$, then $|f| \in \mathsf{C}$; $\textbf{iii)}$ for all $c
\in \R$, the constant function $f \equiv c$ belongs to $\mathsf{C}$.

From \citep{douc:moulines:2005} we adapt the following definitions.

\begin{definition}[Consistency]
A weighted sample $\{(\parti{0:m}{i}{}, \wgt{m}{i})\}_{i = 1}^{M_N}$
on $\Xset^{m+1}$ is said to be \emph{consistent} for the probability
measure $\mu$ and the (proper) set $\mathsf{C} \subseteq
\Lp{1}(\Xset^{m+1}, \mu)$ if, for any $f \in \mathsf{C}$, as $N
\rightarrow \infty$,
\[
\begin{split}
&(\wgtsum{m})^{-1} \sum_{i=1}^{M_N} \wgt{m}{i} f(\parti{0:m}{i}{})
\stackrel{\Prob}{\longrightarrow} \mu f \eqsp,\\
&(\wgtsum{m})^{-1} \max_{1 \leq i \leq M_N} \wgt{m}{i}
\stackrel{\Prob}{\longrightarrow} 0 \eqsp.
\end{split}
\]
\end{definition}

\begin{definition}[Asymptotic normality]
A weighted sample $\{(\parti{0:m}{i}{}, \wgt{m}{i})\}_{i = 1}^{M_N}$
on $\Xset^{m+1}$ is called \emph{asymptotically normal} (abbreviated a.n.) for
$(\mu, \mathsf{A}, \mathsf{W}, \sigma, \gamma, \{ a_N
\}_{N=1}^\infty)$ if, as $N
\rightarrow
\infty$,
\[
\begin{split}
&a_N (\wgtsum{m})^{-1} \sum_{i=1}^{M_N} \wgt{m}{i} [f(\parti{0:m}{i}{})-\mu f]
\stackrel{\mathcal{D}}{\longrightarrow} \mathcal{N}[0, \sigma^2(f)]
\quad \mathrm{for\ any\ }
f \in \mathsf{A} \eqsp,\\
&a_N^2 (\wgtsum{m})^{-1} \sum_{i=1}^{M_N} (\wgt{m}{i})^2 f(\parti{0:m}{i}{})
\stackrel{\Prob}{\longrightarrow} \gamma f \quad \mathrm{for\
  any\ } f \in \mathsf{W} \eqsp,\\
&a_N (\wgtsum{m})^{-1} \max_{1 \leq i \leq M_N} \wgt{m}{i}
\stackrel{\Prob}{\longrightarrow} 0 \eqsp.
\end{split}
\]
\end{definition}
The main contribution of this section is the following results, which
establish consistency and asymptotic normality of weighted samples
produced by the \TSS and \SSAPF algorithms. For all $k \geq 0$, we define a
transformation $\smoothop{k}$ on the set of
$\smooth{k}{k}$-integrable functions by
\begin{equation} \label{eq:def:smoothop}
\smoothop{k}[f](\vect{x}_{0:k}) \define f(\vect{x}_{0:k}) -
\smooth{k}{k} f \eqsp, \quad \vect{x}_{0:k} \in \Xset^{k+1} \eqsp.
\end{equation}
In addition, we impose the following assumptions.

\begin{hyp}{hyp:CLT:assumption}
For all $k \geq 1$, $\fstimpfunc{k} \in
\Lp{2}(\Xset^{k+1}, \smooth{k}{k})$ and $\sdimpfunc{k} \in
\Lp{1}(\Xset^{k+1}, \smooth{k}{k})$, where $\fstimpfunc{k}$ and
$\sdimpfunc{k}$ are defined in \eqref{eq:def:fstimpfunc} and
\eqref{eq:def:sdimpfunc}, respectively.
\end{hyp}

\begin{hyp}{CLT:initial-sample}
\begin{itemize}
\item[{\it i)}]
$\mathsf{A}_0 \subseteq \Lp{1}(\Xset, \smooth{0}{0})$ is a proper set and
$\sigma_0 : \mathsf{A}_0 \rightarrow \R^+$ is a function satisfying,
for all $f \in \mathsf{A}_0$ and $a \in \R$, $\sigma_0(a f) = |a|
\sigma_0(f)$.
\item[{\it ii)}]
 The initial sample $\{(\parti{0}{i}{},1)\}_{i = 1}^N$ is consistent
 for $[ \Lp{1}(\Xset, \smooth{0}{0}),
\smooth{0}{0} ]$ and a.n. for $[ \smooth{0}{0},
  \mathsf{A}_0, \mathsf{W}_0, \sigma_0, \gamma_0, \{ \sqrt{N}
  \}_{N=1}^\infty ]$.
\end{itemize}
\end{hyp}
\begin{theorem} \label{th:CLT:TSS}
Assume \refhyp{hyp:CLT:assumption} and \refhyp{CLT:initial-sample}
with $(\mathsf{W}_0, \gamma_0) = [\Lp{1}(\Xset, \smooth{0}{0}),
\smooth{0}{0}]$. In
the setting of Algorithm~\ref{alg:TSS}, suppose that the limit $\beta
\define \lim_{N \rightarrow \infty} N / M_N$ exists, where
$\beta \in [0, 1]$. Define recursively the family $\{
\mathsf{A}_k \}_{k = 1}^\infty$ by
\begin{multline} \label{eq:def:A:TSS:CLT}
\mathsf{A}_{k+1} \define \Big \{ f \in \Lp{2}(\Xset^{k+2},
\smooth{k+1}{k+1}):
\pathprop{k}(\cdot, \sdimpfunc{k+1}|f|) \uk{k}(\cdot, |f|) \in
\Lp{1}(\Xset^{k+1},
\smooth{k}{k}), \\
\uk{k}(\cdot, |f|) \in
\mathsf{A}_k \cap \Lp{2}(\Xset^{k+1}, \smooth{k}{k}),
\sdimpfunc{k+1} f^2 \in \Lp{1}(\Xset^{k+2},\ \smooth{k+1}{k+1})
\Big \} \eqsp.
\end{multline}
Furthermore, define recursively the family
$\{ \sigma_k \}_{k = 1}^\infty$ of functionals $\sigma_k : \mathsf{A}_k
\rightarrow \R^+$ by
\begin{multline} \label{eq:var:CLT:TSS}
\sigma_{k+1}^2(f) \define \smooth{k+1}{k+1} \smoothop{k+1}^2 [f] +
\frac{\sigma_k^2 \{ 
  \uk{k}(\cdot, \smoothop{k+1}[f]) \} + \beta \smooth{k}{k} \{
   \fstimpfunc{k} \pathprop{k}( \cdot, \sdimpfunc{k+1}^2
   \smoothop{k+1}^2[f] ) \} \, \smooth{k}{k} 
   \fstimpfunc{k}}{[\smooth{k}{k} \uk{k}(\Xset^{k+2})]^2}
  \eqsp.
 \end{multline}
 Then each $\mathsf{A}_k$ is a proper set for all $k \geq 1$. Moreover,
 each sample $\{(\parti{0:k}{i}{},1)\}_{i = 1}^N$ produced by
 Algorithm~\ref{alg:TSS} is consistent for $[\Lp{1}(\Xset^{k+1},
 \smooth{k}{k}), \smooth{k}{k} ]$ and asymptotically normal
 for $[\smooth{k}{k}, \mathsf{A}_k, \Lp{1}(\Xset^{k+1},
 \smooth{k}{k}), \sigma_k, \smooth{k}{k},
 \{ \sqrt{N} \}_{N=1}^\infty]$.
\end{theorem}

The proof is found in Appendix~\ref{section:appendix:A},
and as a by-product a similar result for the SSAPF
(Algorithm~\ref{alg:single:step:aux}) is obtained. 

\begin{theorem} \label{th:CLT:SSAPF}
Assume \refhyp{hyp:CLT:assumption} and
\refhyp{CLT:initial-sample}. Define the families $\{
\tilde{\mathsf{W}}_k \}_{k = 0}^\infty$ and $\{ \tilde{\mathsf{A}}_k \}_{k =
  0}^\infty$ by
\[
\tilde{\mathsf{W}}_k \define \Big\{ f \in \Lp{1}(\Xset^{k+1},
\smooth{k}{k}): \sdimpfunc{k+1} f \in \Lp{1}(\Xset^{k+1},
\smooth{k}{k}) \Big\} \eqsp,  \quad \tilde{\mathsf{W}}_0 \define
\mathsf{W}_0 \eqsp,
\]
and, with $\tilde{\mathsf{A}}_{0} \define \mathsf{A}_0$,
\begin{multline} \label{eq:def:A:SSAPF:CLT}
\tilde{\mathsf{A}}_{k+1} \define \Big \{ f \in \Lp{1}(\Xset^{k+2},
\smooth{k+1}{k+1}):
\pathprop{k}(\cdot, \sdimpfunc{k+1}|f|) \uk{k}(\cdot, |f|) \in
\Lp{1}(\Xset^{k+1}, \smooth{k}{k}), \\
\uk{k}(\cdot, |f|) \in \tilde{\mathsf{A}}_k, [\uk{k}(\cdot, |f|)]^2 \in
\tilde{\mathsf{W}}_k, \sdimpfunc{k+1} f^2 \in \Lp{1}(\Xset^{k+2},\
\smooth{k+1}{k+1}) 
\Big \} \eqsp.
\end{multline}
Furthermore, define recursively the family
$\{ \tilde{\sigma}_k \}_{k = 0}^\infty$ of functionals
$\tilde{\sigma}_k : \mathsf{A}_k \rightarrow \R^+$ by
\begin{equation} \label{eq:var:CLT:SSAPF}
\tilde{\sigma}_{k+1}^2(f) \define \frac{\tilde{\sigma}_k^2 \{
  \uk{k}(\cdot, \smoothop{k+1}[f]) \} + \smooth{k}{k} \{
   \fstimpfunc{k} \pathprop{k}( \cdot, \sdimpfunc{k+1}^2
   \smoothop{k+1}^2[f] ) \} \, \smooth{k}{k} 
   \fstimpfunc{k}}{[\smooth{k}{k} \uk{k}(\Xset^{k+2})]^2}
  \eqsp, \quad \tilde{\sigma}_0 \define \sigma_0 \eqsp,
 \end{equation}
and the measures $\{ \tilde{\gamma}_k \}_{k = 1}^{\infty}$ by
\[
\tilde{\gamma}_{k+1} f \define \frac{\smooth{k+1}{k+1}(\sdimpfunc{k+1} f) \,
  \smooth{k}{k} \fstimpfunc{k}}{\smooth{k}{k}
  \uk{k}(\Xset^{k+2})} \eqsp, \quad f \in \tilde{\mathsf{W}}_{k+1} \eqsp.
\]
 Then each $\tilde{\mathsf{A}}_k$ is a proper set for all $k \geq
 1$. Moreover, each sample $\{(\parttilde{0:k}{i}{}, \wgttilde{k}{i})\}_{i =
   1}^N$ produced by Algorithm~\ref{alg:single:step:aux} is consistent for
 $[\Lp{1}(\Xset^{k+1}, \smooth{k}{k}), \smooth{k}{k} ]$ and
 asymptotically normal
 for $[\smooth{k}{k}, \tilde{\mathsf{A}}_k, \tilde{\mathsf{W}}_k,
 \tilde{\sigma}_k, \tilde{\gamma}_k, \{ \sqrt{N} \}_{N=1}^\infty]$.
\end{theorem}

Under the assumption of bounded likelihood and second-stage
importance weight functions $\lf{k}$ and $\sdimpfunc{k}$, one can show
that the CLTs stated in Theorems~\ref{th:CLT:TSS} and
\ref{th:CLT:SSAPF} indeed include any functions having
finite second moments \wrt the joint smoothing distributions; that is,
under these assumptions the supplementary constraints on the sets
\eqref{eq:def:A:TSS:CLT} and \eqref{eq:def:A:SSAPF:CLT} are
automatically fulfilled. This is the contents of the statement below.
\begin{hyp}{assumption:bdd:likelihood}
For all $k \geq 0$, $\supnm{\Xset}{\lf{k}} < \infty$ and
$\supnm{\Xset^{k+1}}{\sdimpfunc{k}} < \infty$.
\end{hyp}
\begin{corollary}
\label{cor:A:is:L2}
Assume \refhyp{assumption:bdd:likelihood} and let $\{
\mathsf{A}_k \}_{k = 0}^\infty$ and $\{ \tilde{\mathsf{A}}_k \}_{k =
  0}^\infty$ be defined by \eqref{eq:def:A:TSS:CLT} 
and \eqref{eq:def:A:SSAPF:CLT}, respectively, with $\tilde{\mathsf{A}}_0
= \mathsf{A}_0 \define \Lp{2}(\Xset, \smooth{0}{0})$. Then, for all $k
\geq 1$, $\mathsf{A}_k = \Lp{2}(\Xset^{k+1}, \smooth{k}{k})$ and
$\Lp{2}(\Xset^{k+1}, \smooth{k}{k}) \subseteq \tilde{\mathsf{A}}_k$.
\end{corollary}

For a proof, see Section~\ref{section:proof:corollary:A:is:L2}.

Interestingly, the expressions of $\tilde{\sigma}_{k+1}^2(f)$ and
$\sigma_{k+1}^2(f)$ differ, for $\beta = 1$, \emph{only on the additive term}
$\smooth{k+1}{k+1} \smoothop{k+1}^2 [f]$, that is, the variance of $f$ under
$\smooth{k+1}{k+1}$. This quantity represents the cost of
introducing the second-stage resampling pass, which was proposed as a
mean for preventing the particle approximation from degenerating. In
the coming Section~\ref{section:bounds:on:Lp:error:and:bias} we will
however show that the approximations produced
by the \SSAPF are already stable for a finite time horizon, and that
additional resampling is superfluous. Thus, there are indeed reasons
for strongly questioning whether second-stage resampling should be
performed at all, at least when the same number of particles are used
in the two stages.

\subsection{Bounds on $\Lp{p}$ error and bias}
\label{section:bounds:on:Lp:error:and:bias}

In this part we examine, under suitable regularity conditions and
for a finite particle population, the errors of the approximations
obtained by the \APF in terms $\Lp{p}$ bounds and bounds on the
bias. We preface our main result with some definitions and
assumptions. Denote by $\bddonX{m}$ space of bounded measurable
functions on $\Xset^m$ furnished with the supremum norm $\| f
\|_{\Xset^m, \infty} \define \sup_{\vect{x} \in \Xset^m}
|f(\vect{x})|$. Let, for $f \in \bddonX{m}$, the \emph{oscillation
  semi-norm} (alternatively termed the \emph{global modus of
  continuity}) be defined by $\osc(f) \define
\sup_{(\vect{x},\vect{x}') \in \Xset^{m} \times \Xset^{m}}|f(\vect{x})
- f(\vect{x}')|$. Furthermore, the $\Lp{p}$ norm of a stochastic
variable $X$ is denoted by $\| X \|_p \define \E^{1/p}[|X|^p]$. When
considering sums, we will make use of the standard convention
$\sum_{k=a}^b c_k = 0$ if $b < a$.

In the following we will assume that all measures $\hk{}(x,
\cdot)$, $x \in \Xset$, have densities $\hd{}(x, \cdot)$ \wrt a common
dominating measure $\refm$ on ($\Xset$, $\Xsigm$). Moreover, we suppose that
the following holds.
 \begin{hyp}{assumption:regularity:q}
\begin{itemize}
\item[\it{i)}] $\lb \define \inf_{(x,x') \in \Xset^2} \hd{} (x,x') >
  0$, $\ub \define \sup_{(x,x') \in \Xset^2} \hd{} (x,x') < \infty$.
\item[\it{ii)}] For all $y \in \Yset$,  $\int_\Xset \lf{} (y|x) \,
  \refm( \ud x) > 0$.
\end{itemize}
\end{hyp}
Under \refhyp{assumption:regularity:q} we define
\begin{equation} \label{eq:def:rho}
\rho \define 1 - \frac{\lb}{\ub} \eqsp.
\end{equation}
\begin{hyp}{assumption:finite:weights}
For all $k \geq 0$, $\supnm{\Xset^{k+1}}{\fstimpfunc{k}} < \infty$.
\end{hyp}
Assumption~\refhyp{assumption:regularity:q} is now standard and is often
satisfied when the state space $\Xset$ is compact and implies that
the hidden chain, when evolving conditionally on the observations, is
geometrical ergodic with a mixing rate given by $\rho < 1$. For
comprehensive treatments of such stability properties within the
framework of state space models we refer
to \citet{delmoral:2004}. Finally, let $\mathcal{C}_i(\Xset^{n+1})$ be
the set of bounded measurable functions $f$ on $\Xset^{n+1}$ of type
$f(\vect{x}_{0:n}) = \bar{f}(\vect{x}_{i:n})$ for some function
$\bar{f} : \Xset^{n-i+1} \rightarrow \R$. In this setting we have the following
result, which is proved in
Section~\ref{section:proof:main:deviation:theorem}.
\begin{theorem} \label{th:main:deviation:theorem}
Assume \refhyp{assumption:bdd:likelihood}, \refhyp{assumption:regularity:q},
\refhyp{assumption:finite:weights}, and let $f \in
\mathcal{C}_i(\Xset^{n+1})$ for $0 \leq i \leq n$. Let $\{
(\parttilde{0:k}{i}{}, \wgttilde{k}{i}{}) \}_ {i = 1}^{R_N(r)}$ be a
weighted particle sample produced by Algorithm~$r$, $r = \{1 ,2\}$,
with $R_N(r) \define \ind \{ r = 1 \} M_N + \ind \{ r = 2 \} N$. Then
the following holds true for all $N \geq 1$ and $r = \{1 ,2\}$.
\begin{itemize}
\item[\it{ i) }] For all $p \geq 2$,
\begin{multline*}
\left \|  (\tilde{\Omega}_n^N)^{-1} \sum_{j=1}^{R_N(r)} \wgttilde{n}{j}
    f_i(\parttilde{0:n}{j}{}) - \smooth{n}{n}f_i \right \|_p \\
 \leq B_p \frac{\osc(f_i)}{1 - \rho} \left[ \frac{1}{\lb \sqrt{R_N(r)}}
   \sum_{k=1}^n \frac{\left \| \sdimpfunc{k} \right
  \|_{\Xset^{k+1}, \infty} \left \| \fstimpfunc{k-1} \right
  \|_{\Xset^k, \infty}}{\refm \lf{k}} \rho^{0 \vee (i-k)} \right.\\
+ \Bigg. \frac{\ind \{ r = 1 \}}{\sqrt{N}}
 \left( \frac{\rho}{1 - \rho} +
    n-i \right) +  \frac{\left \| \sdimpfunc{0} \right
  \|_{\Xset, \infty}}{\nu \lf{0} \sqrt{N}}\rho^i \Bigg] \eqsp,
\end{multline*}
\item[\it{ ii) }]
\begin{multline*}
\left| \E \left[ (\tilde{\Omega}_n^N)^{-1} \sum_{j=1}^{R_N(r)} \wgttilde{n}{j}
    f_i(\parttilde{0:n}{j}{}) \right] - \smooth{n}{n}f_i \right| \\
 \leq B \frac{\osc(f_i)}{(1 - \rho)^2} \left[ \frac{1}{R_N(r) \lb^2}
   \sum_{k=1}^n \frac{\left \| \sdimpfunc{k} \right
  \|_{\Xset^{k+1}, \infty}^2 \left \| \fstimpfunc{k-1} \right
  \|_{\Xset^k, \infty}^2}{(\refm \lf{k})^2} \rho^{0 \vee (i-k)} \right.\\
+ \Bigg. \frac{\ind \{ r = 1 \}}{N}
\left( \frac{\rho}{1 - \rho} +  n-i \right) + \frac{\left
    \| \sdimpfunc{0} \right
  \|_{\Xset, \infty}^2}{N (\nu \lf{0})^2} \rho^i \Bigg] \eqsp.
\end{multline*}
\end{itemize}
Here $\rho$ is defined in \eqref{eq:def:rho}, and $B_p$ and $B$ are
universal constants such that $B_p$ depends on $p$ only.
\end{theorem}

Especially, applying, under the assumption that all fractions $\|
\sdimpfunc{k} \|_{\Xset^{k+1}, \infty} \| \fstimpfunc{k-1}
\|_{\Xset^k, \infty} / \refm \lf{k}$ are uniformly bounded in $k$,
Theorem~\ref{th:main:deviation:theorem} for $i = n$, yields error
bounds on the approximate filter distribution which are
\emph{uniformly bounded} in $n$. From this it is obvious that the 
first-stage resampling pass is enough to preserve the sample
stability. Indeed, by avoiding second-stage selection according to
Algorithm~\ref{alg:single:step:aux} we can, since
the middle terms in the bounds above cancel in this case, obtain
even \emph{tighter} control of the $\Lp{p}$ error for a fixed number
of particles.

\section{Identifying asymptotically optimal first-stage weights}
\label{section:optimal:first-stage:weights}

The formulas \eqref{eq:var:CLT:TSS} and \eqref{eq:var:CLT:SSAPF} for
the asymptotic variances of the \TSS and \SSAPF may look complicated
at a first sight, but by carefully examining the same we will obtain
important knowledge of how to choose the first-stage importance weight
functions $\fstimpfunc{k}$ in order
to robustify the \APF.

Assume that we have run
the \APF up to time $k$ and are about to design suitable first-stage
weights for the next iteration. In this setting,
we call a first-stage weight function $\fstimpfunc{k}'[f]$,
possibly depending on the target function
$f \in \mathsf{A}_{k+1}$ and satisfying \refhyp{hyp:CLT:assumption},
\emph{optimal} (at time $k$) if it provides a minimal increase
of asymptotic variance at a single iteration of the \APF algorithm, that is, if
$\sigma^2_{k+1} \{ \fstimpfunc{k}'[f] \}(f) \leq \sigma^2_{k+1}\{ t
\}(f)$ (or $\tilde{\sigma}^2_{k+1} \{ \fstimpfunc{k}'[f] \}(f) \leq
\tilde{\sigma}^2_{k+1}\{ t \}(f)$) for all other measurable and
positive weight functions $t$. Here we let $\sigma^2_{k+1}\{ t \}(f)$
denote the asymptotic variance induced by $t$. Define, for
$\vect{x}_{0:k} \in \Xset^{k+1}$, 
\begin{equation} \label{eq:otimal:Tk}
\fstimpfunc{k}^\ast [f](\vect{x}_{0:k}) \define \sqrt{\int_\Xset
 \lf{k+1}^2(x_{k+1}) \left[ \frac{\ud \hk{} (x_k, \cdot)}{\ud \prop{k}(x_k,
  \cdot)}(x_{k+1}) \right]^2
 \smoothop{k+1}^2[f](\vect{x}_{0:k+1}) \, \prop{k}(x_k, \ud x_{k+1})}
  \eqsp,
\end{equation}
and let $\sdimpfunc{k+1}^\ast[f]$ denote the second-stage importance
weight function induced by $\fstimpfunc{k}^\ast [f]$ according to
\eqref{eq:def:sdimpfunc}. We are now ready to state the main result of
this section. The proof is found in
Section~\ref{section:proof:optimal:weights}.
\begin{theorem}
\label{th:optimal:weights}
Let $k \geq 0$ and define $\fstimpfunc{k}^\ast$ by
\eqref{eq:otimal:Tk}. Then the following is valid.
\begin{itemize}
\item[{\it i)}]
Let the assumptions of Theorem~\ref{th:CLT:TSS} hold and suppose that
$f \in \{ f' \in \mathsf{A}_{k+1} : \fstimpfunc{k}^\ast [f'] \in
\Lp{2}(\Xset^{k+1}, \smooth{k}{k}), \sdimpfunc{k+1}^\ast [f'] \in
\Lp{1}(\Xset^{k+2}, \smooth{k+1}{k+1})\}$. Then $\fstimpfunc{k}^\ast$ is
optimal for Algorithm~\ref{alg:TSS} and the corresponding minimal
variance is given by
\begin{equation*} 
\sigma_{k+1}^2 \{ \fstimpfunc{k}^\ast \}(f) = \smooth{k+1}{k+1}
\smoothop{k+1}^2 [f] + \frac{\sigma_k^2[
  \uk{k}(\cdot, \smoothop{k+1}[f])] +
\beta (\smooth{k}{k} \fstimpfunc{k}^\ast[f])^2}{[\smooth{k}{k}
  \uk{k}(\Xset^{k+2})]^2} \eqsp.
\end{equation*}
\item[{\it ii)}]
Let the assumptions of Theorem~\ref{th:CLT:SSAPF} hold and suppose that
$f \in \{ f' \in \tilde{\mathsf{A}}_{k+1} : \fstimpfunc{k}^\ast [f'] \in
\Lp{2}(\Xset^{k+1}, \smooth{k}{k}), \sdimpfunc{k+1}^\ast [f'] \in
\Lp{1}(\Xset^{k+2}, \smooth{k+1}{k+1})\}$. Then $\fstimpfunc{k}^\ast$ is
optimal for Algorithm~\ref{alg:single:step:aux} and the corresponding minimal
variance is given by
\begin{equation*} 
\tilde{\sigma}_{k+1}^2 \{ \fstimpfunc{k}^\ast \}(f) = \frac{\tilde{\sigma}_k^2[
  \uk{k}(\cdot, \smoothop{k+1}[f])] +
(\smooth{k}{k} \fstimpfunc{k}^\ast[f])^2}{[\smooth{k}{k}
  \uk{k}(\Xset^{k+2})]^2} \eqsp.
\end{equation*}
\end{itemize}
\end{theorem}

The functions $\fstimpfunc{k}^\ast$ have a natural interpretation in terms of
optimal sample allocation for \emph{stratified sampling}. Consider the
mixture $\pi = \sum_{i=1}^d w_i \mu_i$, each $\mu_i$ being a measure
on some measurable space $(\mathsf{E}, \mathcal{E})$ and $\sum_{i=1}^d
w_i = 1$, and the problem of estimating, for some given
$\pi$-integrable target function $f$, the expectation $\pi f$. In
order to relate this to the particle filtering paradigm, we will make
use of Algorithm~\ref{alg:biased:stratified:sampling}.
\begin{algorithm}[h]
\caption{Stratified importance sampling}
\label{alg:biased:stratified:sampling}
\begin{algorithmic}[1]
\For{$i = 1, \ldots, N$}
\State draw an index $J_i$ multinomially \wrt $\tau_j$, $1
  \leq j \leq d$, $\sum_{j=1}^d \tau_j = 1$;
\State simulate $\xi_i \sim \nu_{J_i}$, and
\State compute the weights $\omega_i \define \left. \frac{w_j}{\tau_j}
  \frac{\ud \mu_j}{\ud \nu_j} \right|_{j=J_i}$
\EndFor
\State Take $\{(\xi_i, \omega_i)\}_{i=1}^N$ as an approximation of $\pi$.
\end{algorithmic}
\end{algorithm}
In other words, we perform Monte Carlo estimation of $\pi f$ by means of
sampling from some proposal mixture $\sum_{j=1}^d \tau_j
\nu_j$ and forming a self-normalised estimate---cf. the technique applied in
Section~\ref{section:two-stage:sampling} for sampling from
$\smoothmixture{k+1}{k+1}$. In this setting, the following CLT can be
established under weak assumptions: 
\begin{equation*}
\sqrt{N} \left[ \sum_{i = 1}^N \frac{\omega_i}{\sum_{\ell = 1}^N
    \omega_\ell}f(\xi_i) - \pi f \right]
\stackrel{\mathcal{D}}{\longrightarrow}
\mathcal{N} \left[0, \sum_{j=1}^d \frac{w_j^2
    \alpha_j(f)}{\tau_j} \right] \eqsp,
\end{equation*}
with, for $x \in \mathsf{E}$,
\[
\alpha_i(f) \define \int_{\mathsf{E}} \left[ \frac{\ud \mu_i}{\ud
    \nu_i} (x) \right]^2 \Pi^2[f](x) \, \nu_i(\ud x) \quad
\textrm{and} \quad
\Pi[f](x) \define f(x) - \pi f \eqsp.
\]
Minimising the asymptotic variance $\sum_{i=1}^d [w_i^2 \alpha_i(f) /
\tau_i]$ \wrt $\tau_i$, $1 \leq i \leq d$, e.g., by means of the
Lagrange multiplicator method (the details are simple), yields
the optimal weights
\[
\tau_i^\ast \propto w_i \sqrt{\alpha_i (f)} = w_i \sqrt{
 \int_{\mathsf{E}} \left[ \frac{\ud \mu_i}{\ud
    \nu_i} (x) \right]^2 \Pi^2[f](x) \, \nu_i(\ud x)} \eqsp,
\]
and the similarity between this expression and that of the optimal
first-stage importance weight functions $\fstimpfunc{k}^\ast$ is
striking. This strongly supports the idea of interpreting optimal
sample allocation for particle filters in terms of variance reduction
for stratified sampling.


\section{Implementations}
\label{section:implementations}
As shown in the previous section, the utilisation of the optimal
weights \eqref{eq:otimal:Tk} provides, for a given sequence $\{
\prop{k} \}_{k=0}^{\infty}$ of proposal kernels, the most efficient of
all particle filters belonging to the large class covered by
Algorithm~\ref{alg:single:step:aux} (including the standard bootstrap
filter and any fully adapted particle filter). However, exact
computation of the optimal weights is in general infeasible by two
reasons: firstly, they depend (via $\smoothop{k+1}[f]$) on the
expectation $\smooth{k+1}{k+1} f$, that is, the quantity that we aim
to estimate, and, secondly, they involve the evaluation of a
complicated integral. A comprehensive treatment of the important issue
of how to approximate the optimal weights is beyond the scope of this
paper, but in the following three examples we discuss some possible heuristics
for doing this. 

\subsection{Nonlinear Gaussian model}
\label{section:nonlinear:gaussian:model}
In order to form an initial idea of the performance of the optimal
SSAPF in practice, we apply the method to a first order (possibly
nonlinear) autoregressive model observed in noise: 
\begin{equation} \label{eq:nonlinear:gaussian:SSM}
\begin{split}
X_{k+1} & = m(X_k) + \sigma_w(X_k) W_{k+1} \eqsp, \\
Y_k &= X_k + \sigma_v V_k \eqsp,
\end{split}
\end{equation}
with $\{ W_k \}_{k=1}^\infty$ and $\{ V_k \}_{k=0}^\infty$ being
mutually independent sets of standard normal distributed variables
such that $W_{k+1}$ is independent of $(X_i, Y_i)$, $0 \leq i \leq k$,
and $V_k$ is independent of $X_k$, $(X_i, Y_i)$, $0 \leq i \leq
k-1$. Here the functions $\sigma_w : \R \rightarrow \R^+$ and $m : \R
\rightarrow \R$ are measurable, and $\Xset = \R$. As observed by
\citet{pitt:shephard:1999}, it is, for all models of form
\eqref{eq:nonlinear:gaussian:SSM}, possible to propose new particle
using the optimal kernel directly, yielding 
$\pathprop{k} = \opt{k}$ and, for $(x,x') \in \R^2$,
\begin{equation} \label{eq:opt:prop}
\propdens{k}(x,  x') = \frac{1}{\sigmaopt{k}(x) \sqrt{2 \pi}} \exp
\left \{- \frac{[x' - \meanopt{k}(x)]^2}{2 \sigmaopt{k}^2 (x)} \right
\} \eqsp, 
\end{equation}
with $\propdens{k}$ denoting the density of $\prop{k}$ \wrt the
Lebesque measure, and 
\begin{equation} \label{eq:mean:opt:kernel}
\meanopt{k}(x) \define \left[ \frac{y_{k+1}}{\sigma_v^2} +
  \frac{m_k(x)}{\sigma_w^2(x)} \right] \sigmaopt{k}^2(x) \eqsp, \quad
\sigmaopt{k}^2(x) \define \frac{\sigma_v^2 \sigma_w^2(x)}{\sigma_v^2 +
  \sigma_w^2(x)} \eqsp. 
\end{equation}
For the proposal \eqref{eq:opt:prop} it is, for $\vect{x}_{k:k+1} \in
\R^2$, valid that 
\begin{equation} \label{eq:fa:first:stage:weights}
\lf{k+1}(x_{k+1}) \frac{\ud \hk{}(x_k, \cdot)}{\ud \prop{k}(x_k,
  \cdot)}(x_{k+1}) \propto h_k(x_k) \define
\frac{\sigmaopt{k}(x_k)}{\sigma_w(x_k)} \exp \left[
  \frac{\meanopt{k}^2(x_k)}{2 \sigmaopt{k}^2(x_k)} - \frac{m^2(x_k)}{2
    \sigma_w^2(x_k)} \right] \eqsp, 
\end{equation}
and since the right hand side does not depend on $x_{k+1}$ we can, by
letting $\fstimpfunc{k}(\vect{x}_{0:k}) = h_k(x_k)$, $\vect{x}_{0:k}
\in \R^{k+1}$, obtain second-stage weights being indeed unity
(providing a sample of genuinely
$\smoothmixture{k+1}{k+1}$-distributed particles). 
When this is achieved, \citet{pitt:shephard:1999} call the particle
filter \emph{fully adapted}. There is however nothing in the previous
theoretical analysis that supports the idea that aiming at evenly
distributed second-stage weights is always convenient, and this will
also be illustrated in the simulations below. On the other hand, it is
possible to find \emph{cases} when the fully adapted particle filter
is very close to being optimal; see again the following discussion. 

In this part we will study the following two special cases of
\eqref{eq:nonlinear:gaussian:SSM}:

\begin{itemize}
\item{$m(X_k) \equiv \phi X_k$ and $\sigma_w(X_k) \equiv \sigma$}

For a linear/Gaussian model of this kind, exact expressions of the
optimal weights can be obtained using the Kalman filter. We set $\phi
= 0.9$ and let the latent chain be put at stationarity from the
beginning, that is, $X_0 \sim \mathcal{N}[0,\sigma^2/(1 -
\phi^2)]$. In this setting, we simulated, for $\sigma = \sigma_v =
0.1$, a record $\vect{y}_{0:10}$ of observations and estimated the
filter posterior means (corresponding to projection target functions
$\proj{k}(\vect{x}_{0:k}) \define x_k$, $\vect{x}_{0:k} \in \R^{k+1}$) along
this trajectory by applying (1) \SSAPF based on true optimal weights,
(2) \SSAPF based on the generic weights $\PSimpfunct{k}$ of
\citet{pitt:shephard:1999}, and (3) the standard bootstrap particle
filter (that is, \SSAPF with $\fstimpfunc{k} \equiv 1$). In this first
experiment, the prior kernel $\hk{}$ was taken as proposal in all
cases, and since the optimal weights are derived using asymptotic
arguments we used as many as $100,\!000$ particles for all
algorithms. The result is displayed in
Figure~\ref{fig:lin:gauss:aux:vs:boot:a}, and it is clear that
operating with true optimal allocation weights improves---as
expected---the MSE performance in comparison with the other methods. 

The main motivation of \citet{pitt:shephard:1999} for introducing
auxiliary particle filtering was to robustify the particle
approximation to outliers. Thus, we mimic \citet[Example
7.2.3]{cappe:moulines:ryden:2005} and repeat the experiment above for
the observation record $\vect{y}_{0:5} = (-0.652, -0.345, -0.676,
1.142, 0.721, 20)$, standard deviations $\sigma_v = 1$, $\sigma =
0.1$, and the smaller particle sample size $N = 10,\!000$. Note the
large discrepancy of the last observation $y_5$, which in this case is
located at a distance of 20 standard deviations from the mean of the
stationary distribution. The outcome is plotted in
Figure~\ref{fig:lin:gauss:aux:vs:boot:b} from which it is evident that
the particle filter based on the optimal weights is the most efficient
also in this case; moreover, the performance of the standard auxiliary
particle filter is improved in comparison with the bootstrap
filter. Figure~\ref{fig:weight:functions} displays a plot of the
weight functions $\fstimpfunc{4}^\ast$ and $\PSimpfunct{4}$ for the
same observation record. It is clear that $\PSimpfunct{4}$ is not too
far away from the optimal weight function (which is close to symmetric
in this extreme situation) in this case, even if the distance between
the functions as measured with the supremum norm is still
significant. 

\begin{figure}
\centering
\subfigure[]{
\label{fig:lin:gauss:aux:vs:boot:a}
\includegraphics[width=.45\textwidth]{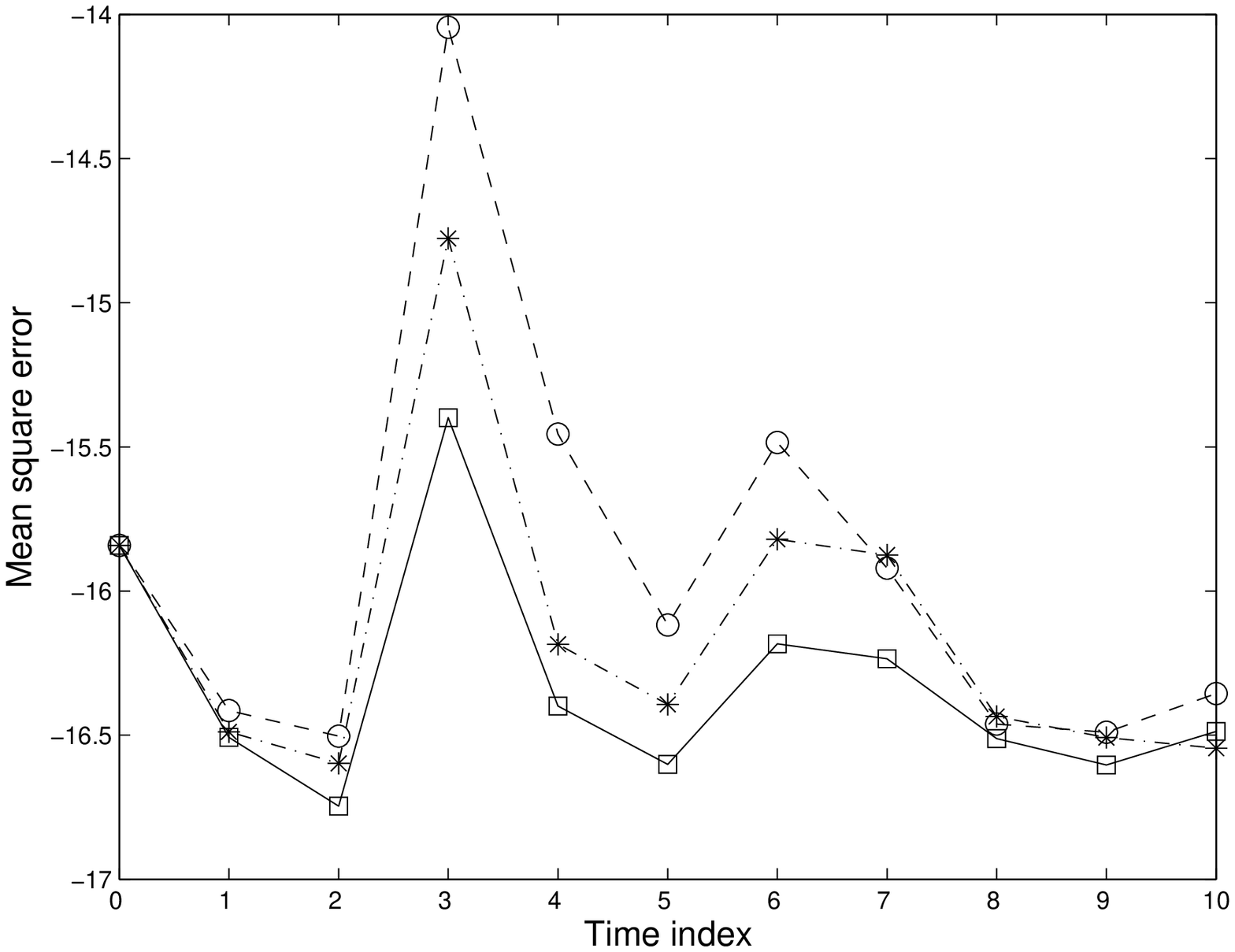}
}
\subfigure[]{
\label{fig:lin:gauss:aux:vs:boot:b}
\includegraphics[width=.45\textwidth]{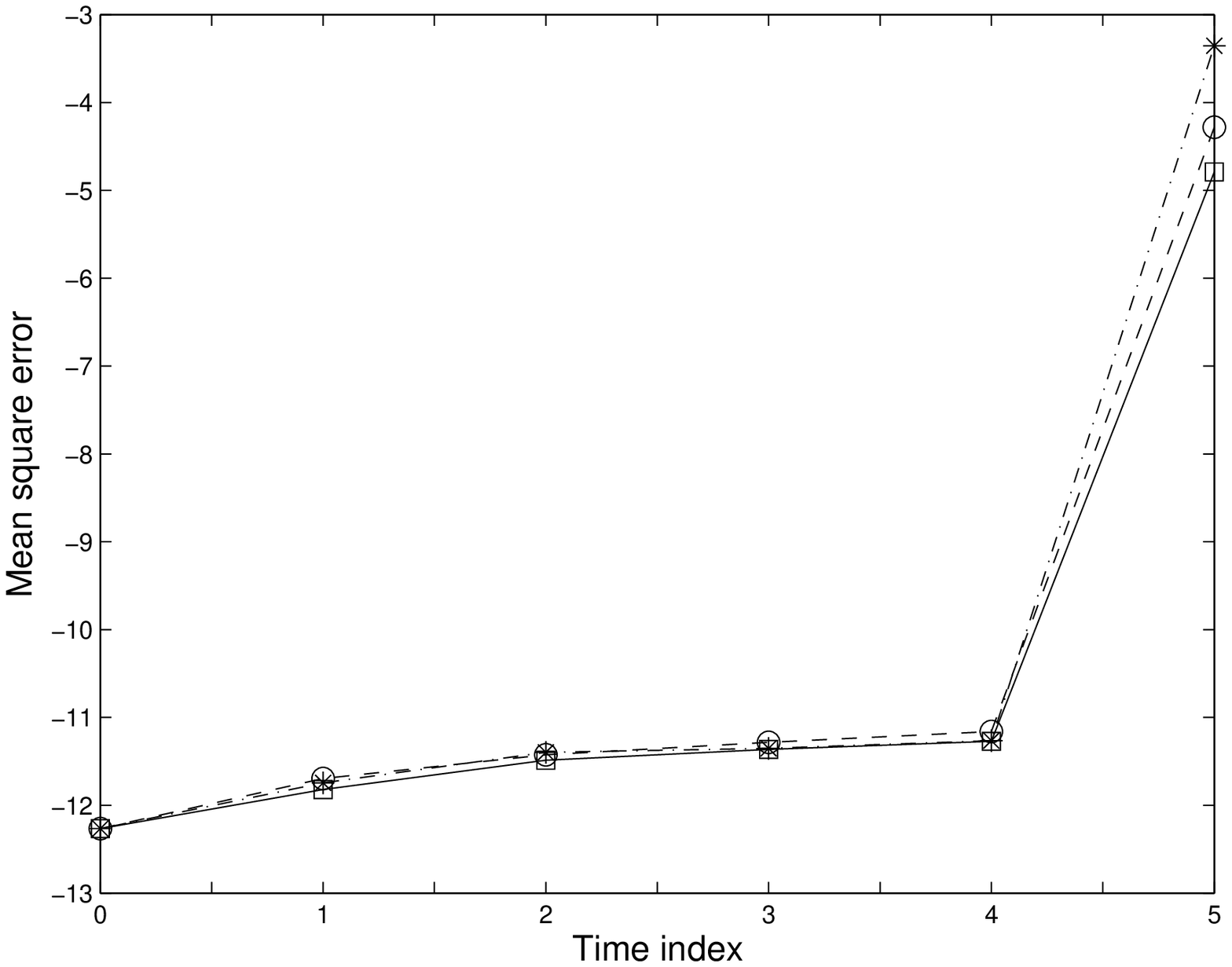}
}
\caption{Plot of MSE perfomances (on log-scale) of the bootstrap
  particle filter ($\ast$), the \SSAPF based on optimal weights
  ($\square$), and the \SSAPF based on the generic weights
  $\PSimpfunct{k}$ of \citet{pitt:shephard:1999} ($\circ$). The MSE values
  are founded on 100,000 particles and 400 runs of each algorithm.} 
\end{figure}

\begin{figure}
\centering
\includegraphics[width=.45\textwidth]{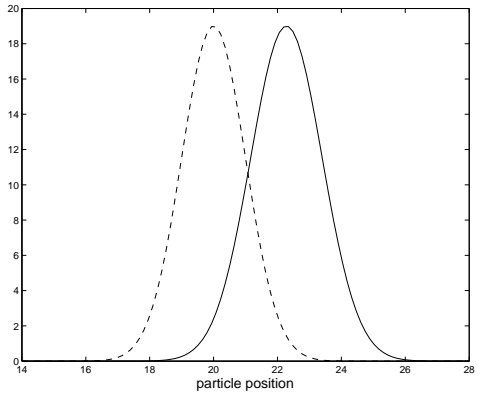}
\caption{Plot of the first-stage importance weight functions
  $\fstimpfunc{4}^\ast$ (unbroken line) and $\PSimpfunct{4}$ (dashed
  line) in the presence of an outlier.} 
\label{fig:weight:functions}
\end{figure}

Finally, we implement the fully adapted filter (with proposal kernels
and first stage-weights given by \eqref{eq:opt:prop} and
\eqref{eq:fa:first:stage:weights}, respectively) and compare this with
the \SSAPF based on the same proposal
\eqref{eq:fa:first:stage:weights} and optimal first-stage weights, the
latter being given by, for $\vect{x}_{0:k} \in \R^{k+1}$ and $h_k$
defined in \eqref{eq:fa:first:stage:weights}, 
\begin{equation} \label{eq:opt:weight:fa:case}
\begin{split}
\fstimpfunc{k}^\ast [\proj{k+1}](\vect{x}_{0:k})
&\propto h_k(x_k) \sqrt{\int_\R
  \smoothop{k+1}^2[\proj{k+1}](\vect{x}_{k+1}) \, \prop{k}(x_k, \ud
  x_{k+1})} \\ 
&= h_k(x_k) \sqrt{\sigmaopt{k}^2(x_k) + \meanopt{k}^2(x_k) - 2
  \meanopt{k}(x_k) \smooth{k+1}{k+1}\proj{k+1} +
  \smooth{k+1}{k+1}^2 \proj{k+1}} 
\end{split}
\end{equation}
in this case. We note that $h_k$, that is, the first-stage weight
function for the fully adapted filter, enters as a factor in the
optimal weight function \eqref{eq:opt:weight:fa:case}. Moreover,
recall the definitions \eqref{eq:mean:opt:kernel} of $\meanopt{k}$ and
$\sigmaopt{k}$; in the case of very informative observations,
corresponding to $\sigma_v \ll \sigma$, it holds that $\sigmaopt{k}(x)
\approx \sigma_v$ and $\meanopt{k}(x) \approx y_{k+1}$ with good
precision for moderate values of $x \in \R$ (that is, values not too
far away from the mean of the stationary distribution of $X$). Thus,
the factor beside $h_k$ in \eqref{eq:opt:weight:fa:case} is more or
less constant in this case, implying that the fully adapted and
optimal first-stage weight filters are close to equivalent. This
observation is perfectly confirmed in
Figure~\ref{fig:fa:vs:opt:filter:a} which presents MSE performances
for $\sigma_v = 0.1$, $\sigma = 1$, and $N = 10,\!000$. In the same
figure, the bootstrap filter and the standard auxiliary filter based
on generic weights are included for a comparison, and these
(particularly the latter) are marred with significantly larger Monte
Carlo errors. On the contrary, in the case of non-informative
observations, that is, $\sigma_v \gg \sigma$, we note that
$\sigmaopt{k}(x) \approx \sigma$, $\meanopt{k}(x) \approx \phi x$ and
conclude that the optimal kernel is close the prior kernel $\hk{}$. In
addition, the exponent of $h_k$ vanishes, implying uniform first-stage
weights for the fully adapted particle filter. Thus, the fully adapted
filter will be close to the bootstrap filter in this case, and
Figure~\ref{fig:fa:vs:opt:filter:b} seems to confirm this
remark. Moreover, the optimal first-stage weight filter does clearly
better than the others in terms of MSE performance. 
\begin{figure}
\centering
\subfigure[]{
\label{fig:fa:vs:opt:filter:a}
\includegraphics[width=.45\textwidth]{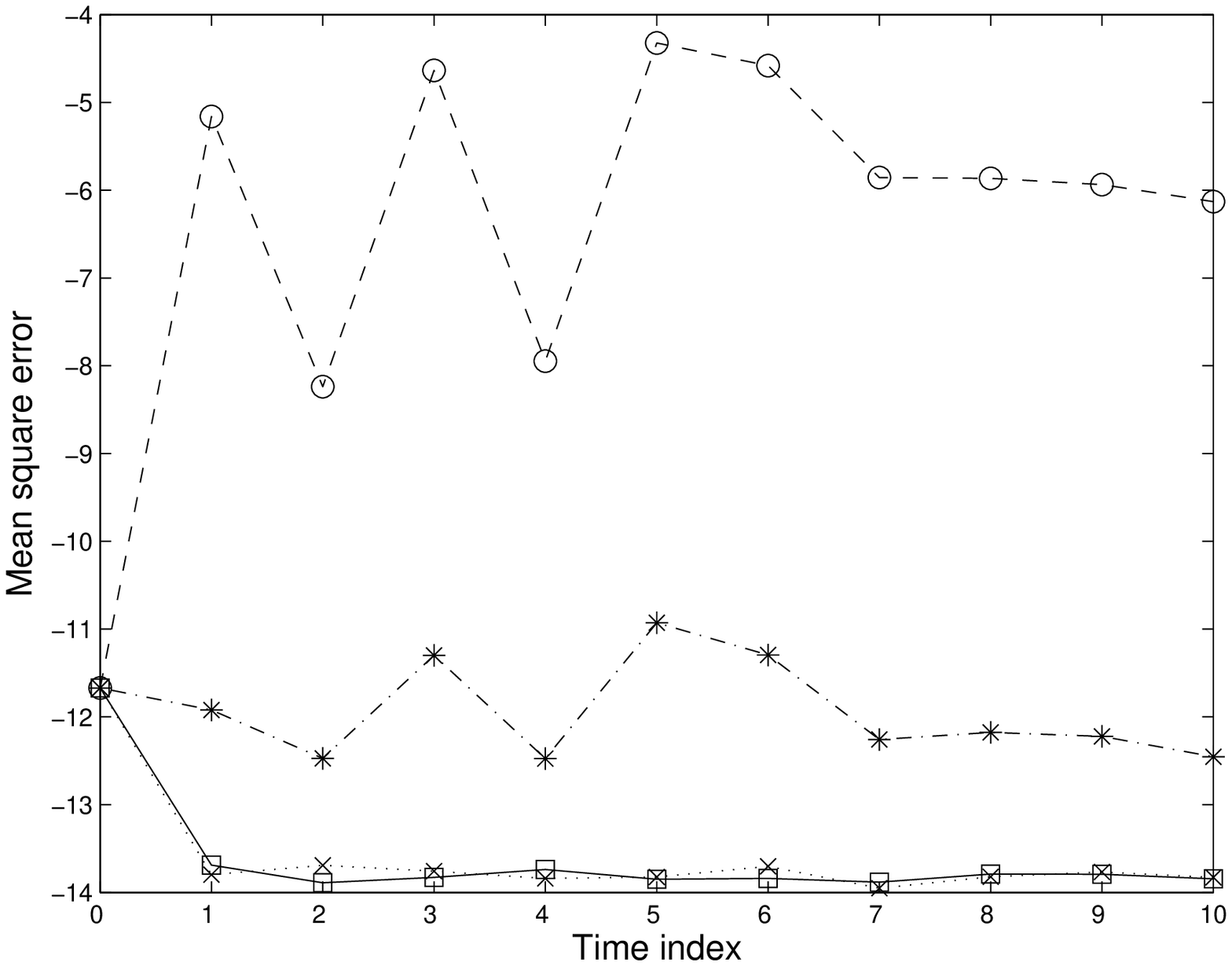}
}
\subfigure[]{
\label{fig:fa:vs:opt:filter:b}
\includegraphics[width=.45\textwidth]{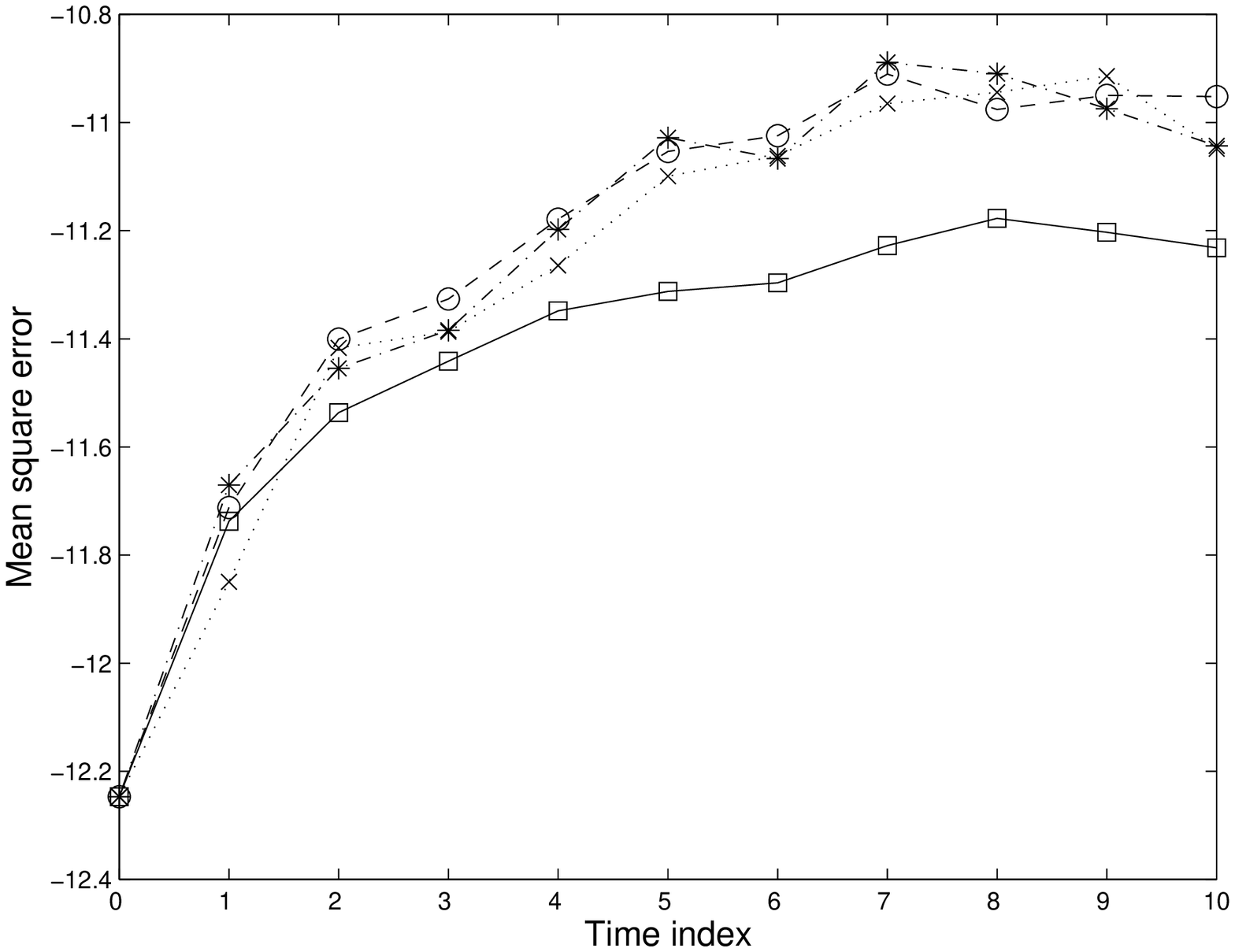}
}
\caption{Plot of MSE perfomances (on log-scale) of the bootstrap
  particle filter ($\ast$), the \SSAPF based on optimal weights
  ($\square$), the \SSAPF based on the generic weights
  $\PSimpfunct{k}$ ($\circ$), and the fully adapted \SSAPF ($\times$) for
  the Linear/Gaussian model in
  Section~\ref{section:nonlinear:gaussian:model}. The MSE values are
  computed using 10,000 particles and 400 runs of each algorithm.} 
\end{figure}


\item{$m(X_k) \equiv 0$ and $\sigma_w(X_k) \equiv \sqrt{\beta_0 +
      \beta_1 X_k^2}$} 

Here we deal with the classical Gaussian autoregressive conditional
heteroscedasticity (ARCH) model
\citep[see][]{bollerslev:engle:nelson:1994} observed in noise. Since
the nonlinear state equation precludes exact computation of the
filtered means, implementing the optimal first-stage weight \SSAPF is
considerably more challenging in this case. The problem can however be
tackled by means of an introductory \emph{zero-stage} simulation pass,
based on $R \ll N$ particles, in which a crude estimate
of $\smooth{k+1}{k+1}f$ is obtained. For instance, this can be
achieved by applying the standard bootstrap filter with multinomial
resampling. Using this approach, we computed again MSE values for the
bootstrap filter, the standard \SSAPF based on generic weights, the
fully adapted SSAPF, and the (approximate) optimal first-stage weight
SSAPF, the latter using the optimal proposal kernel. Each algorithm
used $10,\!000$ particles and the number of particles in the prefatory
pass was set to $R = N/10 = 1000$, implying only a minor additional
computational work. An imitation of the true filter means was obtained
by running the bootstrap filter with as many as $500,\!000$
particles. In compliance with the foregoing, we considered the case of
informative (Figure~\ref{fig:ARCH:fa:vs:boot:a}) as well as
non-informative (Figure~\ref{fig:ARCH:fa:vs:boot:b}) observations,
corresponding to $(\beta_0, \beta_1, \sigma_v) = (9,5,1)$ and
$(\beta_0, \beta_1, \sigma_v) = (0.1,1,3)$, respectively. Since
$\sigmaopt{k}(x) \approx \sigma_v$, $\meanopt{k}(x) \approx y_{k+1}$
in the latter case, we should, in accordance with the previous
discussion, again expect the fully adapted filter to be close to that
based on optimal first-stage weights. This is also confirmed in the
plot. For the former parameter set, the fully adapted \SSAPF exhibits
a MSE performance close to that of the bootstrap filter, while the
optimal first-stage weight SSAPF is clearly superior. 
\begin{figure}
\centering
\subfigure[]{
\label{fig:ARCH:fa:vs:boot:a}
\includegraphics[width=.45\textwidth]{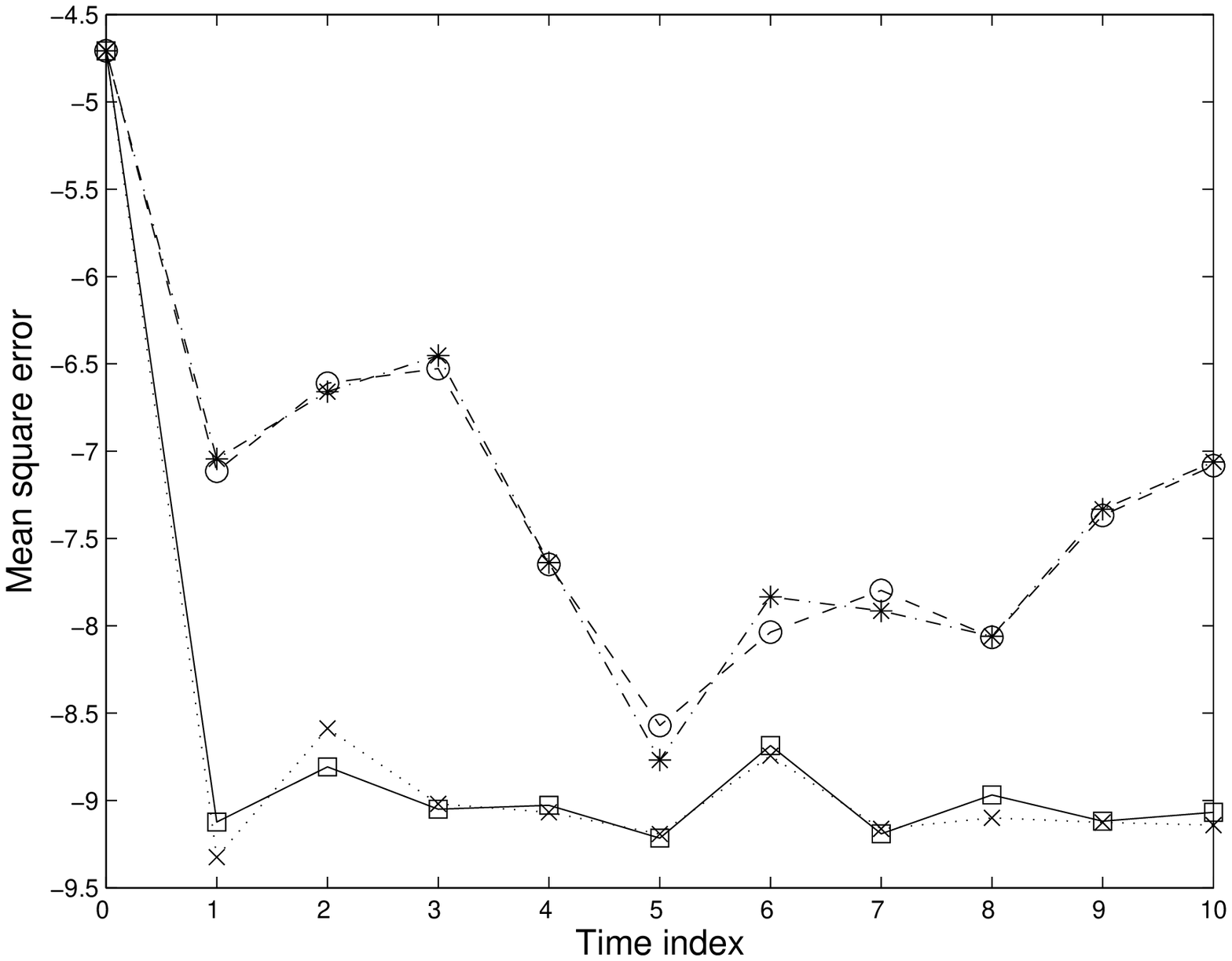}
}
\subfigure[]{
\label{fig:ARCH:fa:vs:boot:b}
\includegraphics[width=.45\textwidth]{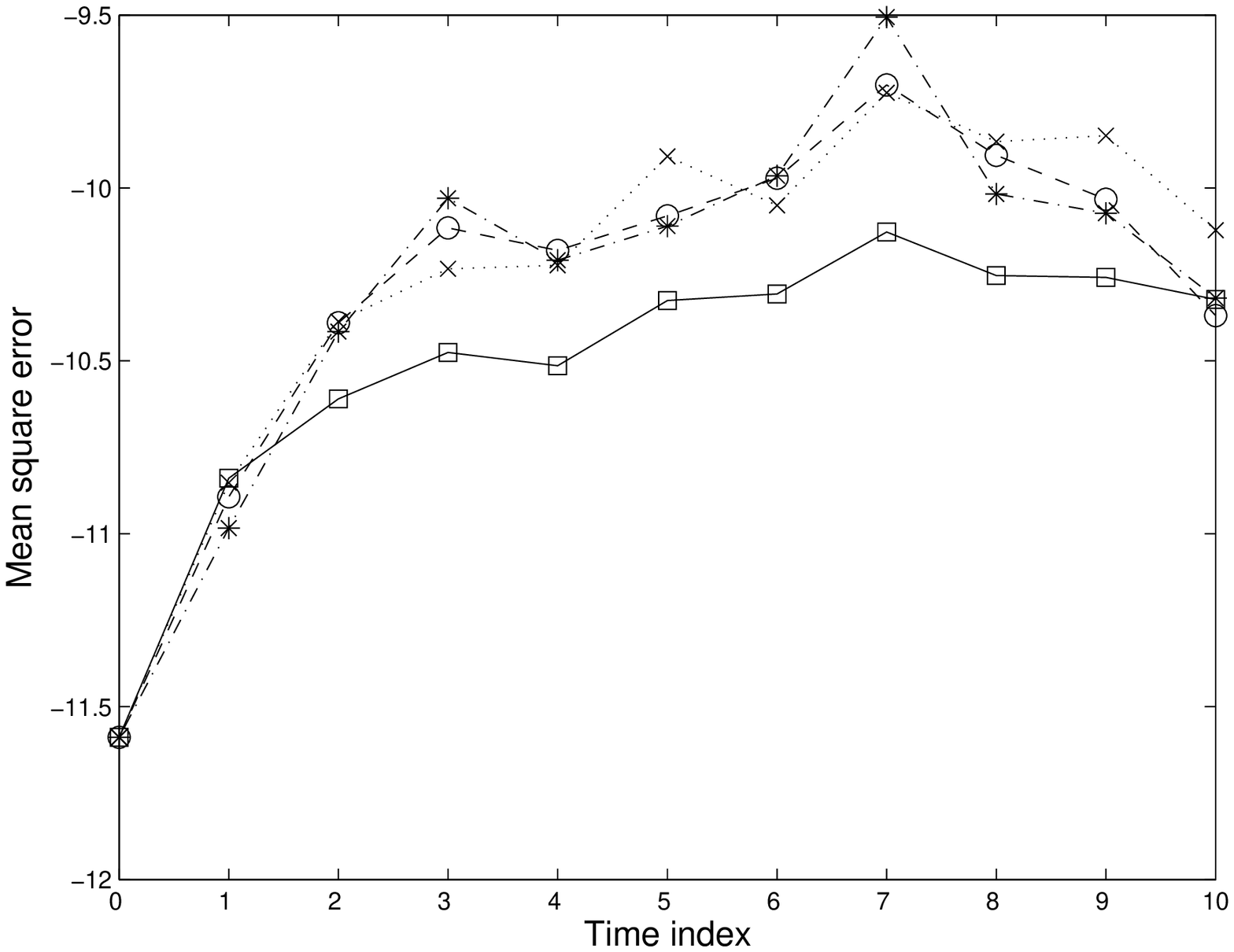}
}
\caption{Plot of MSE perfomances (on log-scale) of the bootstrap
  particle filter ($\ast$), the \SSAPF based on optimal weights
  ($\square$), the \SSAPF based on the generic weights
  $\PSimpfunct{k}$ ($\circ$), and the fully adapted \SSAPF ($\times$) for
  the ARCH model in
  Section~\ref{section:nonlinear:gaussian:model}. The MSE values are
  computed using 10,000 particles and 400 runs of each algorithm.} 
\end{figure}


\end{itemize}

\subsection{Stochastic volatility}
\label{section:SV}
As a final example we consider the canonical discrete-time
\emph{stochastic volatility} (SV) \emph{model} \citep{hull:white:1987}
given by 
\[
\begin{split}
X_{k+1} &= \phi X_k + \sigma W_{k+1} \eqsp, \\
Y_k &= \beta \exp(X_k/2) V_k \eqsp,
\end{split}
\]
where $\Xset = \R$, and $\{ W_k \}_{k=1}^\infty$ and $\{ V_k
\}_{k=0}^\infty$ are as in
Example~\ref{section:nonlinear:gaussian:model}. Here $X$ and $Y$ are
log-volatility and log-returns, respectively, where the former are
assumed to be stationary. Also this model was treated by
\citet{pitt:shephard:1999}, who discussed approximate full adaptation
of the particle filter by means of a second order Taylor approximation
of the concave function $x' \mapsto \log \lf{k+1}(x')$. More
specifically, by multiplying the approximate observation density
obtained in this way with $\hd{}(x,x')$, $(x,x')\in \R^2$, yielding a
Gaussian approximation of the optimal kernel density, nearly even
second-stage weights can be obtained. We proceed in the same spirit,
approximating however directly the (log-concave) function $x' \mapsto
\lf{k+1}(x') \hd{}(x,x')$ by means of a second order Taylor expansion
of $x' \mapsto \log [\lf{k+1}(x') \hd{}(x,x')]$ around the mode
$\meanoptSV{k}(x)$ (obtained using Newton iterations) of the same: 
\[
\lf{k+1}(x') \hd{}(x,x') \approx \propdens{k}^{\mathrm{u}}(x,x')
\define \lf{k+1}[\meanoptSV{k}(x)] \hd{}[x,\meanoptSV{k}(x)] \exp
\left\{ - \frac{1}{2 \sigmaoptSV{k}^2(x)}[x' - \meanoptSV{k}(x)]^2
\right\} \eqsp, 
\]
with \citep[we refer to][pp.~225--228, for
details]{cappe:moulines:ryden:2005} $\sigmaoptSV{k}^2(x)$ being the
inverted negative of the second order derivative, evaluated at
$\meanoptSV{k}(x)$, of $x' \mapsto \log [\lf{k+1}(x')
\hd{}(x,x')]$. Thus, by letting, for $(x,x')\in \R^2$,
$\propdens{k}(x,x') = \propdens{k}^{\mathrm{u}}(x,x')/\int_{\R}
\propdens{k}^{\mathrm{u}}(x,x'') \, \ud x''$, we obtain 
\begin{equation} \label{eq:fa:first:stage:weights:SV}
\lf{k+1}(x_{k+1}) \frac{\ud \hk{}(x_k, \cdot)}{\ud \prop{k}(x_k,
  \cdot)}(x_{k+1}) \approx \int_{\R} \propdens{k}^{\mathrm{u}}(x_k,x')
\, \ud x' \propto \sigmaoptSV{k}(x_k) \lf{k+1}[\meanoptSV{k}(x_k)]
\hd{}[x,\meanoptSV{k}(x_k)] \eqsp, 
\end{equation}
and letting, for $\vect{x}_{0:k} \in \R^{k+1}$,
$\fstimpfunc{k}(\vect{x}_{0:k}) = \sigmaoptSV{k}(x_k)
\lf{k+1}[\meanoptSV{k}(x_k)] \hd{}[x_k,\meanoptSV{k}(x_k)]$ will imply
a nearly fully adapted particle filter. Moreover, by applying the
approximate relation \eqref{eq:fa:first:stage:weights:SV} to the
expression \eqref{eq:otimal:Tk} of the optimal weights, we get
(cf. \eqref{eq:opt:weight:fa:case}) 
\begin{multline}
\label{eq:opt:weight:SV}
\fstimpfunc{k}^\ast [\proj{k+1}](\vect{x}_{0:k}) \approx \int_{\R}
\propdens{k}^{\mathrm{u}}(x_k,x') \, \ud x' \sqrt{ \int_\R
  \smoothop{k+1}^2[\proj{k+1}](\vect{x}) \, \prop{k}(x_k, \ud x)}
\propto \\ \sigmaoptSV{k}(x_k) \lf{k+1}[\meanoptSV{k}(x_k)]
\hd{}[x,\meanoptSV{k}(x_k)] \sqrt{\sigmaoptSV{k}^2(x_k) +
  \meanoptSV{k}^2(x_k) - 2 \meanoptSV{k}(x_k)
  \smooth{k+1}{k+1}\proj{k+1} + \smooth{k+1}{k+1}^2 \proj{k+1}}
\eqsp. 
\end{multline}

\begin{figure}
\centering
\includegraphics[width=.45\textwidth]{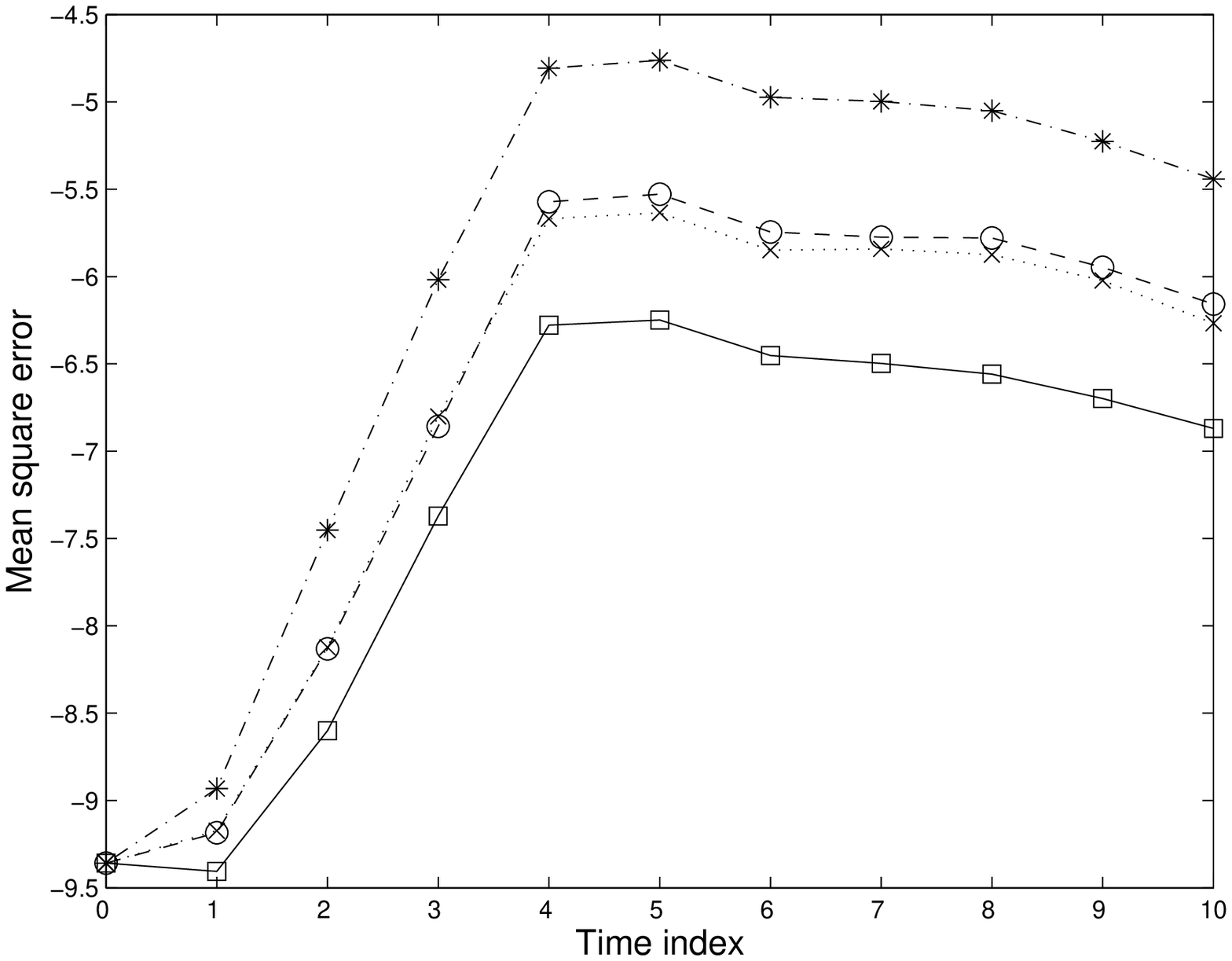}
\caption{Plot of MSE perfomances (on log-scale) of the bootstrap
  particle filter ($\ast$), the \SSAPF based on optimal weights
  ($\square$), the \SSAPF based on the generic weights
  $\PSimpfunct{k}$ ($\circ$), and the fully adapted \SSAPF ($\times$) for
  the SV model in Section~\ref{section:SV}. The MSE values are
  computed using 5,000 particles and 400 runs of each algorithm.} 
\label{fig:SV:fa:vs:boot}
\end{figure}

In this setting, we conducted a numerical experiment where the two
filters above were, again together with the bootstrap filter and the
auxiliary filter based on the generic weights $\PSimpfunct{k}$, run
for the parameters $(\phi, \beta, \sigma) = (0.9702, 0.5992, 0.178)$
\citep[estimated by][from daily returns on the U.~S. dollar against
the U.~K. pound stearling from the first day of trading in 1997 and
for the next 200 days]{pitt:shephard:1999:b}. To make the filtering
problem more challenging, we used a simulated record $\vect{y}_{0:10}$
of observations arising from the initial state $x_0 = 2.19$, being
above the 2\% quantile of the stationary distribution of $X$, implying
a sequence of relatively impetuously fluctuating log-returns. The
number of particles was set to $N = 5,\!000$ for all filters, and the
number of particles used in the prefatory filtering pass (in which a
rough approximation of $\smooth{k+1}{k+1}\proj{k+1}$ in
\eqref{eq:opt:weight:SV} was computed using the bootstrap filter) of
the \SSAPF filter based on optimal first-stage weights was set to $R =
N/5 = 1000$; thus, running the optimal first-stage weight filter is
only marginally more demanding than running the fully adapted
filter. The outcome is displayed in Figure~\ref{fig:SV:fa:vs:boot}. It
is once more obvious that introducing approximate optimal first-stage
weights significantly improves the performance also for the the SV
model, which is recognised as being specially demanding as regards
state estimation. 

\begin{appendix}

\section{Proofs}
\label{section:appendix:A}

\subsection{Proof of Theorem ~\ref{th:CLT:TSS}}
\label{section:proof:theorem:CLT:TSS}

Let us recall the updating scheme described in
Algorithm~\ref{alg:TSS} and formulate it in the following four isolated steps:
\begin{multline} \label{eq:updating:scheme}
\{(\parti{0:k}{i}{},1)\}_ {i = 1}^N \overpil{\textbf{I}: Weighting}
\{(\parti{0:k}{i}{}, \fstwgt{k}{i})\}_{i = 1}^N \overpil{\textbf{II}:
  Resampling (1st stage)} \{(\parthat{0:k}{i}{},1)\}_{i = 1}^{M_N}
\rightarrow \\
\overpil{\textbf{III}: Mutation}
\{(\parttilde{0:k+1}{i}{},\wgttilde{k+1}{i})\}_{i = 1}^{M_N}
\overpil{\textbf{IV}: Resampling (2nd stage)} \{(\parti{0:k+1}{i}{},1)\}_{i =
  1}^N \eqsp,
\end{multline}
where we have set $\parthat{0:k}{i}{} \define
\parti{0:k}{I_k^{N,i}}{}$, $1 \leq i \leq M_N$. Now,
the asymptotic properties stated in Theorem~\ref{th:CLT:TSS} are
established by a chain of applications of
\citep[Theorems~1--4]{douc:moulines:2005}. We will proceed by
induction: assume that the uniformly weighted particle sample
$\{(\parti{0:k}{i}{},1)\}_{i = 1}^N$ is consistent for
$[\Lp{1}(\Xset^{k+1}, \smooth{k}{k}), \smooth{k}{k} ]$ and
asymptotically normal
for $[\smooth{k}{k}, \mathsf{A}_k, \Lp{1}(\Xset^{k+1},
\smooth{k}{k}), \sigma_k, \smooth{k}{k}, \{ \sqrt{N} \}_{N =
  1}^\infty]$, with $\mathsf{A}_k$ being a proper set and $\sigma_k$ such that
$\sigma_k(a f) = |a| \sigma_k(f)$, $f \in \mathsf{A}_k$, $a \in
\R$. We prove, by analysing each of the steps \textbf{(I--IV)}, that
this property is preserved through one iteration of the algorithm.

\textbf{(I)}. Define the measure
\[
\mu_k(A) \define \frac{ \smooth{k}{k}(\fstimpfunc{k}
  \ind_A)}{\smooth{k}{k} \fstimpfunc{k}} \eqsp, \quad A \in
\Xsigm^{\otimes (k+1)} \eqsp.
\]
By applying \citep[Theorem~1]{douc:moulines:2005} for
$R(\vect{x}_{0:k},\cdot) = \delta_{\vect{x}_{0:k}}(\cdot)$,
$L(\vect{x}_{0:k},\cdot) = \fstimpfunc{k}(\vect{x}_{0:k}) \,
\delta_{\vect{x}_{0:k}}(\cdot)$, $\mu = \mu_k$, and $\nu =
\smooth{k}{k}$, we conclude that the weighted sample
$\{(\parti{0:k}{i}{}, \fstwgt{k}{i})\}_{i = 1}^N$ is consistent
for $[\{f \in \Lp{1}(\Xset^{k+1},\mu_k): \fstimpfunc{k} |f| \in
\Lp{1}(\Xset^{k+1},\smooth{k}{k})\}, \mu_k] =
[\Lp{1}(\Xset^{k+1},\mu_k),\mu_k]$. Here the equality is based on the
fact that $\smooth{k}{k}(\fstimpfunc{k} |f|) = \mu_k |f| \,
\smooth{k}{k} \fstimpfunc{k}$, where the second factor on the right
hand side is bounded by Assumption~\refhyp{hyp:CLT:assumption}. In
addition, by applying \citep[Theorem~1]{douc:moulines:2005} we
conclude that $\{(\parti{0:k}{i}{},\fstwgt{k}{i})\}_{i =
  1}^N$ is asymptotically normal for $(\mu_k,
\mathsf{A}_{\mathrm{\mathbf{I}},k}, \mathsf{W}_{\mathrm{\mathbf{I}},k},
\sigma_{\mathrm{\mathbf{I}},k}, \gamma_{\mathrm{\mathbf{I}},k}, \{
\sqrt{N} \}_{N = 1}^\infty)$, where
\[
\begin{split}
\mathsf{A}_{\mathrm{\mathbf{I}},k} & \define \Big \{ f \in
\Lp{1}(\Xset^{k+1},\mu_k):
\fstimpfunc{k} |f| \in \mathsf{A}_k, \fstimpfunc{k} f \in \Lp{2}(\Xset^{k+1},
\smooth{k}{k}) \Big \} \\
& = \Big \{ f \in \Lp{1}(\Xset^{k+1},\mu_k):
\fstimpfunc{k} f \in \mathsf{A}_k \cap \Lp{2}(\Xset^{k+1},
\smooth{k}{k}) \Big
\} \eqsp, \\
\mathsf{W}_{\mathrm{\mathbf{I}},k} &\define \Big \{ f \in
\Lp{1}(\Xset^{k+1},\mu_k): \fstimpfunc{k}^2 |f| \in
\Lp{1}(\Xset^{k+1}, \smooth{k}{k}) \Big \}
\end{split}
\]
are proper sets, and
\[
\begin{split}
\sigma_{\mathrm{\mathbf{I}},k}^2(f) &\define \sigma_k^2 \left[
  \frac{\fstimpfunc{k}(f - \mu_k f)}{\smooth{k}{k} \fstimpfunc{k}}
\right] =  \frac{\sigma_k^2[ \fstimpfunc{k}(f - \mu_k
  f)]}{(\smooth{k}{k} \fstimpfunc{k})^2} \eqsp,
\quad f \in \mathsf{A}_{\mathrm{\mathbf{I}},k} \eqsp,\\
\gamma_{\mathrm{\mathbf{I}},k}f &\define \frac{\smooth{k}{k}(\fstimpfunc{k}^2
  f)}{(\smooth{k}{k} \fstimpfunc{k})^2} \eqsp,
\quad f \in \mathsf{W}_{\mathrm{\mathbf{I}},k} \eqsp.
\end{split}
\]

\textbf{(II)}. By using \citep[Theorems~3 and 4]{douc:moulines:2005} we deduce
that $\{(\parthat{0:k}{i}{},1)\}_{i = 1}^{M_N}$ is
consistent for $[\Lp{1}(\Xset^{k+1},\mu_k), \mu_k]$ and
a.n. 
for $[\mu_k, \mathsf{A}_{\mathrm{\mathbf{II}},k}, \Lp{1}(\Xset^{k+1}, \mu_k),
\sigma_{\mathrm{\mathbf{II}},k}, \beta \mu_{k},$ $\{ \sqrt{N} \}_{N =
  1}^\infty]$, where
\[
\mathsf{A}_{\mathrm{\mathbf{II}},k} \define \Big \{ f \in
\mathsf{A}_{\mathrm{\mathbf{I}},k}:
f \in \Lp{2}(\Xset^{k+1}, \mu_k) \Big \} = \Big \{ f \in
\Lp{2}(\Xset^{k+1}, \mu_k): \fstimpfunc{k} f \in
\mathsf{A}_k \cap \Lp{2}(\Xset^{k+1},
\smooth{k}{k}) \Big \}
\]
is a proper set, and
\[
\sigma^2_{\mathrm{\mathbf{II}},k}(f) \define \beta
\mu_k[(f - \mu_k f)^2] + \sigma^2_{\mathrm{\mathbf{I}},k}(f) = \beta
\mu_k[(f - \mu_k f)^2] 
+ \frac{\sigma_k^2[ \fstimpfunc{k}(f
  - \mu_k f)]}{(\smooth{k}{k} \fstimpfunc{k})^2} \eqsp, \quad f \in
\mathsf{A}_{\mathrm{\mathbf{II}},k} \eqsp.
\]

\textbf{(III)}. 
We argue as in step \textbf{(I)}, but this time for $\nu = \mu_k$, $R
= \pathprop{k}$, and $L(\cdot, A) = \pathprop{k}(\cdot,
\sdimpfunc{k+1} \ind_A)$, $A \in
\Xsigm^{\otimes (k+2)}$, providing the target distribution
\begin{equation} \label{eq:target:dist:III}
\mu (A) = \frac{\mu_k \pathprop{k}(\sdimpfunc{k+1} \ind_A)}{\mu_k
  \pathprop{k} \sdimpfunc{k+1}} = \frac{\smooth{k}{k}
  \uk{k}(A)}{\smooth{k}{k} \uk{k}(\Xset^{k+2})} =
\smooth{k+1}{k+1}(A) \eqsp, \quad A \in \Xsigm^{\otimes (k+2)} \eqsp.
\end{equation}
This yields, applying \citep[Theorems~1 and
2]{douc:moulines:2005}, that
$\{(\parttilde{k+1}{i}{},\wgttilde{k+1}{i})\}_{i = 1}^{M_N}$ is
consistent for
\begin{multline} \label{eq:setid}
\left[ \Big\{ f \in \Lp{1}(\Xset^{k+2}, \smooth{k+1}{k+1}),
\pathprop{k}(\cdot, \sdimpfunc{k+1} |f|) \in \Lp{1}(\Xset^{k+1},
\mu_k) \Big\},\smooth{k+1}{k+1} \right] \\ = \left[ \Lp{1}(\Xset^{k+2},
\smooth{k+1}{k+1}),\smooth{k+1}{k+1} \right] \eqsp,
\end{multline}
where \eqref{eq:setid} follows, since $\mu_k
\pathprop{k}(\sdimpfunc{k+1} |f|) \, \smooth{k}{k} \fstimpfunc{k} =
\smooth{k}{k} \uk{k}(\Xset^{k+2}) \,
\smooth{k+1}{k+1}|f|$, from \refhyp{hyp:CLT:assumption}, and
a.n. for $(\smooth{k+1}{k+1},
\mathsf{A}_{\mathrm{\mathbf{III}},k+1},$
$\mathsf{W}_{\mathrm{\mathbf{III}},k+1},
\sigma_{\mathrm{\mathbf{III}},k+1},
\gamma_{\mathrm{\mathbf{III}},k+1}, \{ \sqrt{N}
\}_{N = 1}^\infty)$. Here
\[
\begin{split}
\lefteqn{\mathsf{A}_{\mathrm{\mathbf{III}},k+1}} \\
&\define  \Big \{ f \in
\Lp{1}(\Xset^{k+2},\smooth{k+1}{k+1}): \pathprop{k}(\cdot,
\sdimpfunc{k+1} |f|) \in \mathsf{A}_{\mathrm{\mathbf{II}},k},
\pathprop{k}(\cdot, \sdimpfunc{k+1}^2 f^2) \in
\Lp{1}(\Xset^{k+1}, \mu_k) \Big \} \\
&= \Big \{ f \in \Lp{1}(\Xset^{k+2}, \smooth{k+1}{k+1}):
\pathprop{k}(\cdot, \sdimpfunc{k+1} |f|) \in \Lp{2}(\Xset^{k+1},
\mu_k),\\
&\hspace{5mm} \fstimpfunc{k} \pathprop{k}(\cdot, \sdimpfunc{k+1}
|f|) \in \mathsf{A}_k \cap \Lp{2}(\Xset^{k+1}, \smooth{k}{k}),
\pathprop{k}(\cdot, \sdimpfunc{k+1}^2 f^2) \in
\Lp{1}(\Xset^{k+1}, \mu_k) \Big \} \\
&=\Big \{ f \in \Lp{1}(\Xset^{k+2}, \smooth{k+1}{k+1}):
\pathprop{k}(\cdot, \sdimpfunc{k+1}|f|) \uk{k}(\cdot, |f|) \in
\Lp{1}(\Xset^{k+1},
\smooth{k}{k}),\\
&\hspace{5mm}
\uk{k}(\cdot, |f|)
\in \mathsf{A}_k \cap \Lp{2}(\Xset^{k+1}, \smooth{k}{k}),
\sdimpfunc{k+1} f^2 \in
\Lp{1}(\Xset^{k+2},\smooth{k+1}{k+1})  \Big\}
\end{split}
\]
and
\[
\begin{split}
\mathsf{W}_{\mathrm{\mathbf{III}},k+1}& \define \Big\{ f \in
\Lp{1}(\Xset^{k+2},
\smooth{k+1}{k+1}): \pathprop{k}(\cdot, \sdimpfunc{k+1}^2|f|) \in
\Lp{1}(\Xset^{k+1},\mu_k) \Big\} \\
&= \Big\{ f \in \Lp{1}(\Xset^{k+2}, \smooth{k+1}{k+1}): \sdimpfunc{k+1}
f \in \Lp{1}(\Xset^{k+2}, \smooth{k+1}{k+1})  \Big\}
\end{split}
\]
are proper sets.
In addition, from the identity \eqref{eq:target:dist:III}
we obtain that
\[
 \mu_k \pathprop{k}( \sdimpfunc{k+1} \smoothop{k+1}[f] ) = 0 \eqsp,
\]
where $\smoothop{k+1}$ is defined in \eqref{eq:def:smoothop}, yielding
\begin{equation*}
\begin{split}
\lefteqn{\sigma_{\mathrm{\mathbf{III}},k+1}^2(f)} \\
& \define \sigma^2_{\mathrm{\mathbf{II}},k}
\left\{ \frac{\pathprop{k} (\cdot,
     \sdimpfunc{k+1}\smoothop{k+1}[f])}{\mu_k \pathprop{k}
     \sdimpfunc{k+1}} \right\} + \frac{\beta \mu_k \pathprop{k}( \{
   \sdimpfunc{k+1}
  \smoothop{k+1}[f] - \pathprop{k}(\cdot, \sdimpfunc{k+1}
  \smoothop{k+1}[f]) \}^2 )}{(\mu_k \pathprop{k} \sdimpfunc{k+1})^2}\\
&= \frac{\beta \mu_k( \{
  \pathprop{k}(\sdimpfunc{k+1}\smoothop{k+1}[f]) \}^2)}{ (\mu_k
  \pathprop{k} \sdimpfunc{k+1})^2} + \frac{\sigma_k^2 \{ \fstimpfunc{k}
  \pathprop{k}(\cdot, \sdimpfunc{k+1}\smoothop{k+1}[f]) \}}{(\smooth{k}{k}
   \fstimpfunc{k})^2(\mu_k \pathprop{k} \sdimpfunc{k+1})^2} \\
&\hspace{30mm} + \frac{\beta \mu_k \pathprop{k}( \{\sdimpfunc{k+1}
  \smoothop{k+1}[f] -  \pathprop{k}(\cdot, \sdimpfunc{k+1}
  \smoothop{k+1}[f]) \}^2 )}{(\mu_k \pathprop{k}
  \sdimpfunc{k+1})^2} \eqsp, \quad f \in \mathsf{A}_{\mathrm{\mathbf{III}},k+1}
\eqsp.
\end{split}
\end{equation*}
Now, applying the equality
\begin{multline*}
 \{ \pathprop{k}( \cdot,
  \sdimpfunc{k+1} \smoothop{k+1}[f]) \}^2 + \pathprop{k}(\cdot,
  \{\sdimpfunc{k+1}\smoothop{k+1}[f] -
   \pathprop{k}(\cdot, \sdimpfunc{k+1} \smoothop{k+1}[f]) \}^2 ) \\
= \pathprop{k}(\cdot, \sdimpfunc{k+1}^2  \smoothop{k+1}^2[f])  \eqsp,
\end{multline*}
provides the variance
\begin{equation} \label{eq:mutation:var}
\sigma_{\mathrm{\mathbf{III}},k+1}^2(f) = \frac{\beta \smooth{k}{k}
  \{ \fstimpfunc{k} \pathprop{k}(\cdot, \sdimpfunc{k+1}^2
  \smoothop{k+1}^2[f]) \} \, \smooth{k}{k}
   \fstimpfunc{k} + \sigma_k^2 \{ \uk{k}(\cdot, \smoothop{k+1}[f])
   \}}{[\smooth{k}{k} \uk{k}(\Xset^{k+2})]^2} \eqsp, \quad f \in
 \mathsf{A}_{\mathrm{\mathbf{III}},k+1}
\eqsp.
\end{equation}

Finally, for $f \in \mathsf{W}_{\mathrm{\mathbf{III}},k+1}$,
\[
\gamma_{\mathrm{\mathbf{III}},k+1}f \define \frac{\beta \mu_k
  \pathprop{k}(\sdimpfunc{k+1}^2 f)}{(\mu_k \pathprop{k}
  \sdimpfunc{k+1})^2} = \frac{\beta \smooth{k+1}{k+1}(\sdimpfunc{k+1} f)
  \, \smooth{k}{k} \fstimpfunc{k}}{\smooth{k}{k}
  \uk{k}(\Xset^{k+2})} \eqsp.
\]

\textbf{(IV)}. The consistency for $[\Lp{1}(\Xset^{k+2},
\smooth{k+1}{k+1}),\smooth{k+1}{k+1}]$ of the uniformly weighted
particle sample $\{(\parti{0:k+1}{i}{},1)\}_{i = 1}^N$ follows from
\citep[Theorem~3]{douc:moulines:2005}. In addition,
applying \citep[Theorem~4]{douc:moulines:2005} yields that the same
sample is a.n. for
$[\smooth{k+1}{k+1}, \mathsf{A}_{\mathrm{\mathbf{IV}},k+1},$
$\Lp{1}(\Xset^{k+2},\smooth{k+1}{k+1}),\sigma_{\mathrm{\mathbf{IV}},k+1},
\smooth{k+1}{k+1}, \{ \sqrt{N} \}_{N = 1}^\infty]$, with
\[
\begin{split}
\mathsf{A}_{\mathrm{\mathbf{IV}},k+1} &\define \Big \{f \in
\mathsf{A}_{\mathrm{\mathbf{III}},k+1}: f \in
\Lp{2}(\Xset^{k+2},\smooth{k+1}{k+1}) \Big\}\\
&= \Big\{ f \in \Lp{2}(\Xset^{k+2}, \smooth{k+1}{k+1}):
\pathprop{k}(\cdot, \sdimpfunc{k+1}|f|) \uk{k}(\cdot, |f|) \in
\Lp{1}(\Xset^{k+1},
\smooth{k}{k}), \\
&\hspace{5mm}
\uk{k}(\cdot, |f|) \in
\mathsf{A}_k \cap \Lp{2}(\Xset^{k+1}, \smooth{k}{k}),
\sdimpfunc{k+1} f^2 \in \Lp{1}(\Xset^{k+2},\smooth{k+1}{k+1})
\Big\}
\end{split}
\]
being proper set, and, for $f \in \mathsf{A}_{\mathrm{\mathbf{IV}},k+1}$,
\[
\begin{split}
\sigma^2_{\mathrm{\mathbf{IV}},k+1}(f) &\define
\smooth{k+1}{k+1} \smoothop{k+1}^2 [f] +
\sigma^2_{\mathrm{\mathbf{III}},k+1}(f) \eqsp, \\
\end{split}
\]
with $\sigma^2_{\mathrm{\mathbf{III}},k+1}(f)$ being defined by
\eqref{eq:mutation:var}. This concludes the proof of the theorem.

\subsection{Proof of Corollary ~\ref{cor:A:is:L2}}
\label{section:proof:corollary:A:is:L2}
We pick $f \in \Lp{2}(\Xset^{k+2},
\smooth{k+1}{k+1})$ and prove that the constraints of the set
$\mathsf{A}_{k+1}$ defined in \eqref{eq:def:A:TSS:CLT} are satisfied
under Assumption~\refhyp{assumption:bdd:likelihood}. Firstly, by
Jensen's inequality,
\[
\begin{split}
\lefteqn{\smooth{k}{k} [\pathprop{k}(\cdot, \sdimpfunc{k+1}|f|)
  \uk{k}(\cdot, |f|)]} \hspace{20mm}\\
&= \smooth{k}{k} \{ \fstimpfunc{k} [\pathprop{k}(\cdot,
  \sdimpfunc{k+1}|f|)]^2 \} \\
&\leq \smooth{k}{k} [\fstimpfunc{k} \pathprop{k}(\cdot,
  \sdimpfunc{k+1}^2 f^2)] \\
&=  \smooth{k}{k} \uk{k} (\sdimpfunc{k+1} f^2) \\
&\leq \supnm{\Xset^{k+2}}{\sdimpfunc{k+1}} \smooth{k}{k} \uk{k}
(\Xset^{k+2}) \, \smooth{k+1}{k+1} (f^2) < \infty \eqsp,
\end{split}
\]
and similarly,
\[
\smooth{k}{k} \{ [ \uk{k}(\cdot, |f|) ]^2 \} \leq
\supnm{\Xset}{\lf{k+1}} \smooth{k}{k} \uk{k}
(\Xset^{k+2}) \, \smooth{k+1}{k+1} (f^2) < \infty \eqsp.
\]
From this, together with the bound
\[
\smooth{k+1}{k+1}(\sdimpfunc{k+1} f^2) \leq
\supnm{\Xset^{k+2}}{\sdimpfunc{k+1}} \smooth{k+1}{k+1} (f^2) < \infty \eqsp,
\]
we conclude that $\mathsf{A}_{k+1} = \Lp{2}(\Xset^{k+2},
\smooth{k+1}{k+1})$.

To prove $\Lp{2}(\Xset^{k+1}, \smooth{k}{k}) \subseteq
\tilde{\mathsf{A}}_k$, note that assumption \refhyp{assumption:bdd:likelihood}
implies $\tilde{\mathsf{W}}_k = \Lp{1}(\Xset^{k+1}, \smooth{k}{k})$
and repeat the arguments above.

\subsection{Proof of Theorem~\ref{th:main:deviation:theorem}}
\label{section:proof:main:deviation:theorem}
Define, for $r \in \{1 ,2\}$ and $R_N(r)$ as defined in Theorem~\ref{th:main:deviation:theorem}, the particle measures
\[
\partsmooth{k}{k}(A) \define \frac{1}{N} \sum_{i=1}^N
  \delta_{\parti{0:k}{i}{}} \quad \text{and} \quad
  \partsmoothtilde{k}{k}(A) \define
\frac{1}{\wgtsumtilde{k}} \sum_{i=1}^{R_N(r)} \wgttilde{k}{i}
\delta_{\parttilde{0:k}{i}{}}(A) \eqsp, \quad A \in \Xsigm^{\otimes
  (k+1)} \eqsp,
\]
playing the role of approximations of the smoothing distribution
$\smooth{k}{k}$. Let $\filt{0}
\define \sigma (\parti{0}{i}{}; 1 \leq i \leq N)$; then the particle
history up to the different steps of loop $m+1$, $m \geq 0$, of
Algorithm~$r$, $r \in \{1,2\}$, is modeled by the filtrations
$\filthat{m} \define \filt{m} \vee \sigma[I_m^{N,i}; 1 \leq i \leq
R_N(r)]$, $\filttilde{m+1} \define
\filt{m} \vee \sigma [\parttilde{0:m+1}{i}{}; 1 \leq i \leq R_N(r)]$, and
\[
\filt{m+1} \define
\begin{cases}
\filttilde{m+1} \vee \sigma (J_{m+1}^{N,i}; 1 \leq i \leq N) \eqsp, &\quad \mathrm{for} \ r = 1 \eqsp, \\
\filttilde{m+1} \eqsp, &\quad \mathrm{for } \ r = 2 \eqsp.
\end{cases}
\]
 respectively.
In the coming proof we will describe one iteration of the \APF
algorithm by the following two operations.

\begin{multline*}
\{(\parti{0:k}{i}{},\wgt{k}{i})\}_ {i = 1}^N
{\raisebox{1.5ex}{
    $\scriptsize \underrightarrow{\hspace{1mm} \text{Sampling from }
      \fstimp{k+1}
      \hspace{1mm}}$}}
\{(\parttilde{0:k+1}{i}{},\wgttilde{k+1}{i})\}_{i = 1}^{R_N(r)} \rightarrow \\
{\raisebox{1.5ex}{
    $\scriptsize \underrightarrow{\hspace{1mm} r = 1\text{: Sampling from }
      \partsmoothtilde{0:k+1}{k+1}
      \hspace{1mm}}$}}
\{(\parti{0:k+1}{i}{},1)\}_{i = 1}^N \eqsp,
\end{multline*}
where, for $A \in \Xsigm^{\otimes (k+2)}$,
\begin{multline} \label{eq:deriv:fstimp}
\fstimp{k+1}(A) \define \Prob \left( \left. \parttilde{0:k+1}{i_0}{} \in A
  \right| \filt{k} \right) = \sum_{j=1}^N \frac{\wgt{k}{j} \fstwgt{k}{j}}{\sum_{\ell=1}^N \wgt{k}{\ell} \fstwgt{k}{\ell}} \, \pathprop{k}(\parti{0:k}{j}{}, A) = \frac{\partsmooth{k}{k}[\fstimpfunc{k}
  \pathprop{k}(\cdot, A )]}{\partsmooth{k}{k}
  \fstimpfunc{k}} \eqsp,
\end{multline}
for some index $i_0 \in \{1, \ldots, R_N(r) \}$ (given $\filt{k}$, the
particles $\parttilde{0:k+1}{i}{}$, $1 \leq i \leq
R_N(r)$, are \iid). Here the initial weights $\{ \wgt{k}{i} \}_{i = 1}^N$ are all equal to one for $r = 1$. The second
operation is valid since, for any $i_0 \in \{1, \ldots, N \}$,
\[
\Prob \left( \left. \parti{0:k+1}{i_0}{} \in A \right| \filttilde{k+1}
\right) = \sum_{j=1}^{R_N(r)} \frac{\wgttilde{k+1}{j}}{\wgtsumtilde{k+1}}
\, \delta_{\parttilde{0:k+1}{j}{}}(A) = \partsmoothtilde{0:k+1}{k+1}(A)
\eqsp, \quad A \in \Xsigm^{\otimes (k+2)} \eqsp.
\]
The fact that the evolution of the particles can be described by two Monte Carlo operations involving conditionally \iid\ variables makes it possible to analyse the error using the
Marcinkiewicz-Zygmund inequality \citep[see][p.~62]{petrov:1995}.

Using this, set, for $1 \leq k \leq n$,
\begin{equation} \label{eq:fsttarg:def}
\fsttarg{k}(A) \define \int_A \frac{\ud \fsttarg{k}}{\ud
  \fstimp{k}}(\vect{x}_{0:k})  \, \fstimp{k}(\ud \vect{x}_{0:k})
\eqsp, \quad A \in \Xsigm^{\otimes (k+1)} \eqsp,
\end{equation}
with, for $\vect{x}_{0:k} \in \Xset^{k+1}$,
\[
\frac{\ud \fsttarg{k}}{\ud
  \fstimp{k}}(\vect{x}_{0:k}) \define \frac{
  \sdimpfunc{k}(\vect{x}_{0:k}) \uk{k} \cdots \uk{n-1}
  (\vect{x}_{0:k}, \Xset^{n+1}) \, \partsmooth{k-1}{k-1}
  \fstimpfunc{k-1}}{\partsmooth{k-1}{k-1} \uk{k-1}
  \cdots \uk{n-1} (\Xset^{n+1})} \eqsp.
\]
Here we apply the standard convention $\uk{\ell} \cdots \uk{m} \define
\operatorname{Id}$ if $m < \ell$. For $k = 0$ we define
\[
\alpha_0(A) \define \int_A \frac{\ud \alpha_0}{\ud
  \varsigma}(x_0) \, \varsigma (\ud x_0) \eqsp,  \quad A \in \Xsigm \eqsp,
\]
with, for $x_0 \in \Xset$,
\[
\frac{\ud \alpha_0}{\ud
  \varsigma}(x_0) \define \frac{\sdimpfunc{0}(x_0) \uk{0} \cdots
\uk{n-1}(x_0, \Xset^{n+1})}{\nu [ \lf{0} \uk{0} \cdots
\uk{n-1}(\cdot, \Xset^{n+1})]} \eqsp.
\]

Similarly, put, for $0 \leq k \leq n-1$,
\begin{equation} \label{eq:sdtarg:def}
\sdtarg{k}(A) \define \int_A \frac{\ud \sdtarg{k}}{\ud
  \partsmoothtilde{k}{k}}(\vect{x}_{0:k})  \,
\partsmoothtilde{k}{k}(\ud \vect{x}_{0:k})
\eqsp, \quad A \in \Xsigm^{\otimes (k+1)} \eqsp,
\end{equation}
where, for $\vect{x}_{0:k} \in \Xset^{k+1}$,
\[
\frac{\ud \sdtarg{k}}{\ud
  \partsmoothtilde{k}{k}}(\vect{x}_{0:k}) \define \frac{\uk{k}
  \cdots \uk{n-1}
  (\vect{x}_{0:k}, \Xset^{n+1})}{\partsmoothtilde{k}{k} \uk{k}
  \cdots \uk{n-1} (\Xset^{n+1})} \eqsp.
\]

The following powerful decomposition is an adaption of a similar one derived by \citet[Lemma~7.2]{cappe:douc:moulines:olsson:2006} (the standard SISR case), being in turn a refinement of a decomposition originally presented by \citet{delmoral:2004}.

\begin{lemma} \label{lemma:key:decomposition}
Let $n \geq 0$. Then, for all $f \in \bddonX{n+1}$, $N \geq 1$, and $r \in \{1, 2\}$,
\begin{equation} \label{eq:key:decomposition}
\partsmoothtilde{0:n}{n}f - \smooth{n}{n}f = \sum_{k=1}^n
A_k^N(f) + \ind \{r = 1\} \sum_{k=0}^{n-1} B_k^N(f) + C^N(f) \eqsp,
\end{equation}
where
 \[
 \begin{split}
 A_k^N (f) & \define \frac{\sum_{i=1}^{R_N(r)} \frac{\ud \fsttarg{k}}{\ud
     \fstimp{k}}(\parttilde{0:k}{i}{})
   \fop{k}{n}[f](\parttilde{0:k}{i}{})}{\sum_{j=1}^{R_N(r)} \frac{\ud
     \fsttarg{k}}{\ud \fstimp{k}}(\parttilde{0:k}{j}{})} - \fsttarg{k}
 \fop{k}{n}[f]
 \eqsp, \\
  B_k^N (f) & \define \frac{\sum_{i=1}^N \frac{\ud \sdtarg{k}}{\ud
        \partsmoothtilde{k}{k}}(\parti{0:k}{i}{})
    \fop{k}{n}[f](\parti{0:k}{i}{})}{\sum_{j=1}^N \frac{\ud
      \sdtarg{k}}{\ud \partsmoothtilde{k}{k}}(\parti{0:k}{j}{})} -
  \sdtarg{k} \fop{k}{n}[f] \eqsp, \\
  C^N (f) &\define \frac{\sum_{i=1}^N \frac{\ud \beta_{0|n}}{\ud
      \varsigma}(\parti{0}{i}{})
    \fop{0}{n}[f](\parti{0}{i}{})}{\sum_{j=1}^N \frac{\ud
      \beta_0}{\ud \varsigma}(\parti{0}{i}{})} - \smooth{n}{n}
  \fop{0}{n}[f] \eqsp,
\end{split}
\]
and the operators $\fop{k}{n} : \bddonX{n+1} \rightarrow
\bddonX{n+1}$, $0 \leq k \leq n$, are, for some fixed points
$\hat{\vect{x}}_{0:k} \in \Xset^{k+1}$, defined by
\begin{equation*} \label{eq:f:hat}
\fop{k}{n}[f] : \vect{x}_{0:k} \mapsto \frac{\uk{k} \cdots
  \uk{n-1} f(\vect{x}_{0:k})}{\uk{k} \cdots \uk{n-1}(\vect{x}_{0:k},
  \Xset^{n+1})} - \frac{\uk{k} \cdots \uk{n-1}
  f(\hat{\vect{x}}_{0:k})}{\uk{k} \cdots
  \uk{n-1}(\hat{\vect{x}}_{0:k}, \Xset^{n+1})} \eqsp.
\end{equation*}
\end{lemma}
\begin{proof}[Proof of Lemma~\ref{lemma:key:decomposition}]
Consider the decomposition
\begin{multline*} \label{eq:key:decomp}
\partsmoothtilde{0:n}{n}f - \smooth{n}{n}f = \sum_{k=1}^n \left[ \frac{\partsmoothtilde{k}{k} \uk{k} \cdots \uk{n-1}
    f}{\partsmoothtilde{k}{k} \uk{k} \cdots \uk{n-1} (\Xset^{n+1})} -
  \frac{\partsmooth{k-1}{k-1} \uk{k-1} \cdots \uk{n-1}
    f}{\partsmooth{k-1}{k-1} \uk{k-1} \cdots \uk{n-1} (\Xset^{n+1})}
\right] \\
 + \ind \{r = 1\} \sum_{k=0}^{n-1} \left[
  \frac{\partsmooth{k}{k} \uk{k} \cdots \uk{n-1}
    f}{\partsmooth{k}{k} \uk{k} \cdots \uk{n-1} (\Xset^{n+1})}
- \frac{\partsmoothtilde{k}{k} \uk{k} \cdots \uk{n-1}
    f}{\partsmoothtilde{k}{k} \uk{k} \cdots \uk{n-1} (\Xset^{n+1})}
\right] \\
+ \frac{\partsmoothtilde{0}{0} \uk{0} \cdots \uk{n-1}
    f}{\partsmoothtilde{0}{0} \uk{0} \cdots \uk{n-1} (\Xset^{n+1})} -
  \smooth{n}{n}f \eqsp.
\end{multline*}
We will show that the three parts of this decomposition are identical with the three parts of \eqref{eq:key:decomposition}. For $k \geq 1$ it holds that, using the definitions \eqref{eq:deriv:fstimp} and \eqref{eq:fsttarg:def} of $\fstimp{k}$ and $\fsttarg{k}$, respectively, and following the lines of \citet[Lemma~7.2]{cappe:douc:moulines:olsson:2006},
\[
\begin{split}
\lefteqn{\frac{\partsmooth{k-1}{k-1} \uk{k-1}
    \cdots \uk{n-1} \uk{n-1} f}{\partsmooth{k-1}{k-1} \uk{k-1} \cdots
    \uk{n-1} (\Xset^{n+1})}} \\
& = \fstimp{k} \left[ \frac{
\sdimpfunc{k}(\cdot) \uk{k}
    \cdots \uk{n-1} f(\cdot) (\partsmooth{k-1}{k-1}
  \fstimpfunc{k-1})}{\partsmooth{k-1}{k-1}
    \uk{k-1} \cdots \uk{n-1} (\Xset^{n+1})} \right] \\
& = \fstimp{k} \left[ \frac{
    \sdimpfunc{k}(\cdot) \uk{k} \cdots
  \uk{n-1} (\cdot, \Xset^{n+1}) (\partsmooth{k-1}{k-1} \fstimpfunc{k-1})}{\partsmooth{k-1}{k-1} \uk{k-1} \cdots
  \uk{n-1}(\Xset^{n+1})} \left\{ \fop{k}{n}[f](\cdot) + \frac{\uk{k} \cdots \uk{n-1}f (\hat{\vect{x}}_{0:k})}{\uk{k} \cdots
  \uk{n-1}(\hat{\vect{x}}_{0:k}, \Xset^{n+1})} \right\}
\right] \\
& = \fsttarg{k} \left[ \fop{k}{n}[f](\cdot) + \frac{\uk{k} \cdots \uk{n-1}f (\hat{\vect{x}}_{0:k})}{\uk{k} \cdots
  \uk{n-1}(\hat{\vect{x}}_{0:k}, \Xset^{n+1})} \right] \\
& = \fsttarg{k} \fop{k}{n}[f] + \frac{\uk{k} \cdots \uk{n-1}f (\hat{\vect{x}}_{0:k})}{\uk{k} \cdots \uk{n-1}(\hat{\vect{x}}_{0:k}, \Xset^{n+1})} \eqsp.
\end{split}
\]
Moreover, by definition,
\[
\frac{\partsmoothtilde{k}{k} \uk{k} \cdots \uk{n-1} f}{\partsmoothtilde{k}{k} \uk{k} \cdots \uk{n-1} (\Xset^{n+1})} = \frac{\sum_{i=1}^{R_N(r)} \frac{\ud \fsttarg{k}}{\ud
     \fstimp{k}}(\parttilde{0:k}{i}{})
   \fop{k}{n}[f](\parttilde{0:k}{i}{})}{\sum_{j=1}^{R_N(r)} \frac{\ud
     \fsttarg{k}}{\ud \fstimp{k}}(\parttilde{0:k}{j}{})} + \frac{\uk{k} \cdots \uk{n-1} f(\hat{\vect{x}}_{0:k})}{\uk{k} \cdots
  \uk{n-1}(\hat{\vect{x}}_{0:k}, \Xset^{n+1})} \eqsp,
\]
yielding
\[
\frac{\partsmoothtilde{k}{k} \uk{k} \cdots \uk{n-1}
    f}{\partsmoothtilde{k}{k} \uk{k} \cdots \uk{n-1} (\Xset^{n+1})} -
  \frac{\partsmooth{k-1}{k-1} \uk{k-1} \cdots \uk{n-1}
    f}{\partsmooth{k-1}{k-1} \uk{k-1} \cdots \uk{n-1} (\Xset^{n+1})} \equiv A_k^N(f) \eqsp.
\]

Similarly, for $r = 1$, using the definition \eqref{eq:sdtarg:def} of $\sdtarg{k}$,
\[
\begin{split}
\frac{\partsmoothtilde{0:k}{k-1} \uk{k-1}
    \cdots \uk{n-1} f}{\partsmoothtilde{0:k}{k-1} \uk{k-1} \cdots
    \uk{n-1} (\Xset^{n+1})}&= \sdtarg{k} \left[ \frac{\uk{k}
    \cdots \uk{n-1} f(\cdot)}{\uk{k} \cdots \uk{n-1} (\Xset^{n+1})} \right] \\
&= \sdtarg{k} \left[ \fop{k}{n}[f](\cdot) + \frac{\uk{k} \cdots \uk{n-1}f (\hat{\vect{x}}_{0:k})}{\uk{k} \cdots
  \uk{n-1}(\hat{\vect{x}}_{0:k}, \Xset^{n+1})} \right] \\
&= \sdtarg{k} \fop{k}{n}[f] + \frac{\uk{k} \cdots \uk{n-1}f (\hat{\vect{x}}_{0:k})}{\uk{k} \cdots \uk{n-1}(\hat{\vect{x}}_{0:k}, \Xset^{n+1})} \eqsp,
\end{split}
\]
and applying the obvious relation
\[
 \frac{\partsmooth{k}{k} \uk{k} \cdots \uk{n-1}
    f}{\partsmooth{k}{k} \uk{k} \cdots \uk{n-1} (\Xset^{n+1})} = \frac{\sum_{i=1}^N \frac{\ud \sdtarg{k}}{\ud
        \partsmoothtilde{k}{k}}(\parti{0:k}{i}{})
    \fop{k}{n}[f](\parti{0:k}{i}{})}{\sum_{j=1}^N \frac{\ud
      \sdtarg{k}}{\ud \partsmoothtilde{k}{k}}(\parti{0:k}{j}{})} + \frac{\uk{k} \cdots \uk{n-1} f(\hat{\vect{x}}_{0:k})}{\uk{k} \cdots
  \uk{n-1}(\hat{\vect{x}}_{0:k}, \Xset^{n+1})} \eqsp,
\]
we obtain the identity
\[
 \frac{\partsmooth{k}{k} \uk{k} \cdots \uk{n-1}
    f}{\partsmooth{k}{k} \uk{k} \cdots \uk{n-1} (\Xset^{n+1})}
- \frac{\partsmoothtilde{k}{k} \uk{k} \cdots \uk{n-1}
    f}{\partsmoothtilde{k}{k} \uk{k} \cdots \uk{n-1} (\Xset^{n+1})} \equiv B_k^N(f) \eqsp.
\]
The equality
\[
 \frac{\partsmoothtilde{0}{0} \uk{0} \cdots \uk{n-1}
    f}{\partsmoothtilde{0}{0} \uk{0} \cdots \uk{n-1} (\Xset^{n+1})} -
  \smooth{n}{n}f \equiv C^N(f)
\]
follows analogously. This completes the proof of the lemma.
\end{proof}

\begin{proof}[Proof of Theorem~\ref{th:main:deviation:theorem}]
From here the proof is a straightforward extension of \citep[][Proposition~7.1]{cappe:douc:moulines:olsson:2006}. To establish part (i), observe that:
\begin{itemize}
\item  A trivial adaption of \citep[][Lemmas~7.3 and~7.4]{cappe:douc:moulines:olsson:2006} gives that
\begin{equation} \label{eq:fop:and:weight:bounds}
\supnm{\Xset^{k+1}}{\fop{k}{n}[f_i]} \leq \osc(f_i) \rho^{0 \vee
  (i-k)} \eqsp, \quad \supnm{\Xset^{k+1}}{\frac{\ud \fsttarg{k}}{\ud
    \fstimp{k}}} \leq \frac{\left \| \sdimpfunc{k} \right
  \|_{\Xset^{k+1}, \infty} \left \| \fstimpfunc{k-1} \right
  \|_{\Xset^{k}, \infty}}{\refm \lf{k} (1-\rho) \lb} \eqsp.
\end{equation}
\item By mimicking the proof of \citep[][Proposition~7.1(i)]{cappe:douc:moulines:olsson:2006}, that is, applying the identity $a/b - c = (a/b)(1-b) + a-c$ to each $A_k^N(f_i)$ and using twice the Marcinkiewicz-Zygmund inequality together with the bounds \eqref{eq:fop:and:weight:bounds}, we obtain the bound
\[
\sqrt{R_N(r)} \left\| A_k^N(f_i) \right\|_p \leq 
B_p \frac{ \osc(f_i) \left \| \sdimpfunc{k} \right
  \|_{\Xset^{k+1}, \infty} \left \| \fstimpfunc{k-1} \right
  \|_{\Xset^{k}, \infty}}{ \refm \lf{k} (1-\rho) \lb}
\rho^{0 \vee (i-k)} \eqsp,
\]
where $B_p$ is a constant depending on $p$ only. We refer to \citep[][Proposition~7.1]{cappe:douc:moulines:olsson:2006} for details.
\item For $r = 1$, inspecting the proof of \citep[][Lemma~7.4]{cappe:douc:moulines:olsson:2006} yields immediately
\[
\supnm{\Xset^{k+1}}{\frac{\ud \sdtarg{k}}{\ud
      \partsmoothtilde{k}{k}}} \leq \frac{1}{1-\rho} \eqsp,
\]
and repeating the arguments of the previous item for $B_k^N(f_i)$ yields
\[
\sqrt{N} \left\| B_k^N(f_i) \right\|_p \leq
B_p \frac{\osc(f_i)}{1-\rho}
\rho^{0 \vee (i-k)} \eqsp.
\]
\item The arguments above apply directly to $C^N(f_i)$, providing

\[
\sqrt{N} \left \| C^N(f_i) \right \|_p \leq B_p \frac{\osc(f_i) \left
    \| \sdimpfunc{0} \right \|_{\Xset, \infty}}{\nu \lf{0}(1-\rho)}
\rho^i \eqsp.
\]
\end{itemize}
We conclude the proof of (i) by summing up.

The proof of (ii) (which mimics the proof of \citep[][Proposition~7.1(ii)]{cappe:douc:moulines:olsson:2006}) follows analogous lines; indeed, repeating the arguments of (i) above for the decomposition $a/b-c = (a/b)(1-b)^2 + (a-c)(1-b) + c(1-b) + a-c$ gives us the bounds
\[
\begin{split}
R_N(r) \left| \E \left[ A_k^N(f_i) \right] \right|& \leq B \frac{ \osc(f_i) \left \| \sdimpfunc{k} \right
  \|_{\Xset^{k+1}, \infty}^2 \left \| \fstimpfunc{k-1} \right
  \|_{\Xset^{k}, \infty}^2}{ (\refm \lf{k})^2 (1-\rho)^2 \lb^2} \rho^{0 \vee
  (i-k)} \eqsp, \\
N \left| \E \left[ B_k^N(f_i) \right] \right| &\leq
B \frac{\osc(f_i)}{(1-\rho)^2}
\rho^{0 \vee (i-k)} \eqsp, \\
 N \left| \E \left[ C^N(f_i) \right]
\right| &\leq B \frac{\osc(f_i) \left
    \| \sdimpfunc{0} \right \|_{\Xset, \infty}^2}{(\nu \lf{0})^2 (1-\rho)^2}
\rho^i \eqsp.
\end{split}
\]
We again refer to \citep[][Proposition~7.1(ii)]{cappe:douc:moulines:olsson:2006} for details, and summing up concludes the proof.
\end{proof}

\subsection{Proof of Theorem~\ref{th:optimal:weights}}
\label{section:proof:optimal:weights}
The statement is a direct implication of H\"older's inequality. Indeed, let
$\fstimpfunc{k}$ be any first-stage importance weight function and write
\begin{equation} \label{eq:optimal:weight:holder:bound}
\begin{split}
(\smooth{k}{k} \fstimpfunc{k}^\ast [f])^2 &= \{ \smooth{k}{k}(
\fstimpfunc{k}^{1/2} \fstimpfunc{k}^{-1/2} \fstimpfunc{k}^\ast[f]) \}^2 \\
&\leq \smooth{k}{k} \fstimpfunc{k} \, \smooth{k}{k} \{
\fstimpfunc{k}^{-1}( \fstimpfunc{k}^\ast [f])^2 \} \eqsp.
\end{split}
\end{equation}
Now the result follows by the formula \eqref{eq:var:CLT:TSS}, the identity
\[
\smooth{k}{k}\{ \fstimpfunc{k}^{-1}( \fstimpfunc{k}^\ast
  [f] )^2 \} = \smooth{k}{k} \{
  \fstimpfunc{k} \pathprop{k}( \cdot, \sdimpfunc{k+1}^2
  \smoothop{k+1}^2[f]) \} \eqsp, \]
and the fact that we have equality in
\eqref{eq:optimal:weight:holder:bound} for $\fstimpfunc{k} = \fstimpfunc{k}^\ast [f]$.

\setcounter{equation}{0}
\end{appendix}
\ \\[3mm]
\textbf{Acknowledgment.} The authors are grateful to Olivier Capp\'{e}
who provided sensible comments on our results that improved the presentation of the paper.




\end{document}